\documentclass[final]{elsarticle}
\makeatletter
\def\ps@pprintTitle{%
 \let\@oddhead\@empty
 \let\@evenhead\@empty
 \def\@oddfoot{}%
 \let\@evenfoot\@oddfoot}
\makeatother
\usepackage{hyperref}
\usepackage[dvipsnames]{xcolor}
\usepackage{amsmath}
\usepackage{amsbsy}
\usepackage{amsthm}
\usepackage{amssymb}
\usepackage[tmargin=1in,bmargin=1in]{geometry}
\usepackage{verbatim}
\usepackage{graphicx}
\usepackage{latexsym}
\usepackage{rotating}
\usepackage{color}
\usepackage{caption}
\usepackage{url}
 \global\bibsep=0pt
\usepackage[normalem]{ulem}
\usepackage{framed,multirow}

\usepackage{url}
\usepackage{xcolor} 
\usepackage[normalem]{ulem}
\usepackage[final]{pdfpages}

\usepackage{algorithm}
\usepackage{float}
\usepackage{siunitx}

\hypersetup{
    colorlinks=true,      		
    linkcolor=blue,       		
    citecolor=blue,       		
    filecolor=black,      		
    urlcolor=red,       		
    bookmarks=false,
    pdffitwindow=true,
    pdfpagelayout=SinglePage
}

\newcommand{\vectorr}[1]{\mathbf{{#1}}}
\newcommand{\npo}{^{n+1}}
\newcommand{\operator}[1]{\mathcal{#1}}

\newcommand{\ie}{{\it i.e. }}
\newcommand{\eg}{{\it e.g. }}
\newcommand{\etal}{{\it et al. }}

\newcommand{\indomain}{\qquad \forall \vectorr{x}\in\Omega^-}

\newcommand{\norm}[1]{\ensuremath{\left\|#1\right\|}}

\newcommand{\ghostinterp}{\operator{I}_{{N}}^{\tilde{N}}}
\newcommand{\proje}{\operator{P}_{\tilde{E}}}
\newcommand{\qproje}{\operator{Q}_{\tilde{E}}}

\newcommand{\ghostinterpedges}{\operator{I}_{{E}}^{\tilde{E}}}
\newcommand{\vertiii}[1]{{\left\vert\kern-0.25ex\left\vert\kern-0.25ex\left\vert #1 
    \right\vert\kern-0.25ex\right\vert\kern-0.25ex\right\vert}}

\definecolor{newcolor}{rgb}{.8,.349,.1}

\theoremstyle{plain}

\begin{document}

 \begin{frontmatter}
      \title{Stable nodal projection method on octree grids}
 \journal{Journal of Computational Physics}

\address[UCMERCED_MATH]{Department of Applied Mathematics, University of California, Merced, California 95343, USA.}
 \author[UCMERCED_MATH]{Matthew Blomquist}
 \author[UCMERCED_MATH]{Scott R. West} 
 \author[UCMERCED_MATH]{Adam L. Binswanger}
 \author[UCMERCED_MATH]{Maxime Theillard} 

\cortext[cor]{Corresponding author: mtheillard@ucmerced.edu}

\begin{abstract}

We propose a novel collocated projection method for solving the incompressible Navier-Stokes equations with arbitrary boundaries. Our approach employs non-graded octree grids, where all variables are stored at the nodes. To discretize the viscosity and projection steps, we utilize supra-convergent finite difference approximations with sharp boundary treatments. We demonstrate the stability of our projection on uniform grids, identify a sufficient stability condition on adaptive grids, and validate these findings numerically. We further demonstrate the accuracy and capabilities of our solver with several canonical two- and three-dimensional simulations of incompressible fluid flows. Overall, our method is second-order accurate, allows for dynamic grid adaptivity with arbitrary geometries, and reduces the overhead in code development through data collocation.

\end{abstract}


\end{frontmatter}

\section{Introduction}

Incompressible fluid flows are ubiquitous in science and engineering applications and lie at the heart of numerous research questions. Developing control strategies to minimize drag \cite{fan2020reinforcement}, understanding arterial wall deformation in the human heart \cite{hsu2014fsiheart}, and even optimizing the energy consumption of automotive spray painting operations \cite{saye2023insights} all require a detailed understanding of incompressible flows. For the majority of these phenomena, analytical solutions do not exist, and experimental approaches can be difficult and costly to create. Numerical simulations are the natural choice for studying these problems. Still, despite decades of computational advancement, their development remains challenging, especially when irregular geometries, adaptive grids, or complex boundary conditions are involved. For this reason, it is essential to develop computational fluid dynamics tools that are accessible and straightforward to implement, which can be achieved, for example, by minimizing the number of unique data structures and simplifying data access patterns. 

The first step in developing a numerical simulation for incompressible flows is typically to discretize the fluid domain using a numerical mesh. This mesh can be represented as a structured set of elements or nodes, as in the Cartesian \cite{harlow1965mac} and curvilinear \cite{KWOK1984curvilinear} styles, or by an unstructured mesh \cite{PAZNER2018344}. Structured meshes are usually easy to generate and can sometimes yield additional accuracy (\eg supra-convergence for nodal finite differences), but special care is needed to handle irregular geometries. Unstructured meshes, conversely, can tessellate irregular geometries, but their generation can be computationally expensive \cite{mullowney2021preparing}. With problems that involve moving geometries or in the context of adaptive grids, unstructured grids often require that computational costs be paid at each time step. For that reason, it is typically preferable to use structured meshes, such as non-graded quad/octrees, and develop the tools to handle irregular geometries. This is the approach we take herein. 

\begin{figure}[b!]
  \centering
    \includegraphics[width = .95\textwidth]{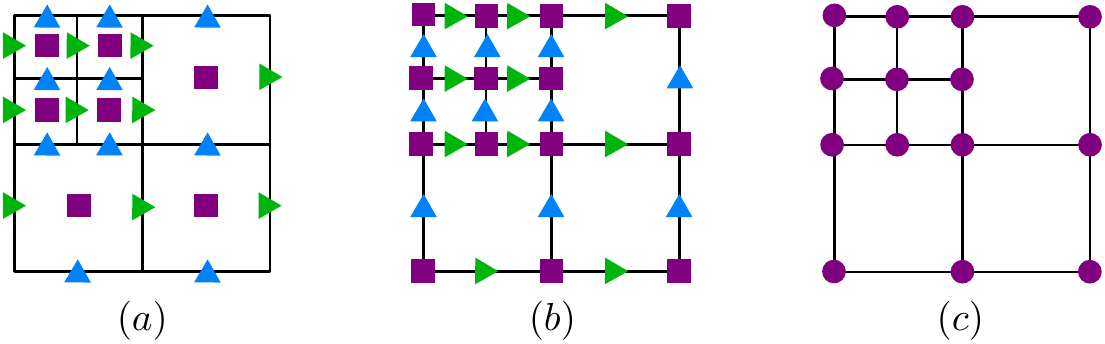}    \caption{Data layouts - (a) Within the Marker and Cell representation the pressure components (\includegraphics[height=0.015\textwidth]{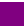}) are stored at the cell centers, the x-velocity (\includegraphics[height=0.015\textwidth]{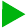}) are stored at the vertical faces, and the y-velocity components (\includegraphics[height=0.015\textwidth]{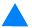}) are stored at the horizontal faces. (b) In \cite{GOMEZ2019478}, Gomez \etal proposed a different staggered storage, where the pressure components are now stored at the nodes, and the velocity components are stored along the edges. (c) In our approach,  all quantities (\includegraphics[height=0.015\textwidth]{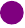}) are collocated at the nodes.    
  }
  \label{fig:mac_gomez_us}
\end{figure}

The primary challenge with numerically solving the incompressible Navier-Stokes equations is related to the coupling between the mass and momentum equations, which manifests as an incompressibility constraint on the velocity field that must be satisfied at each time step. One method for enforcing this constraint is to solve a monolithic system of coupled equations for the velocity and pressure unknowns \cite{ryzhakov2014two}. This approach can often be beneficial when additional physical constraints are coupled with an incompressible flow solver, such as in two-phase flows \cite{aanjaneya2013monolithic,CHO2021110587} and fluid-structure interaction problems \cite{gibou2012monolithic, hyde2019monolithic, CHO2022111304,takahashi2020monolith}. However, the monolithic systems may not be diagonally dominant and thus require the use of expensive solvers (\eg GMRES). For this reason, most computational frameworks for solving the incompressible Navier-Stokes equations use splitting methods.

The projection method, pioneered by Chorin \cite{CHORIN196712} and later proposed as a fractional step method by Temam \cite{temam1969approximation,temam1969approximationII}, decouples the momentum equation and the incompressibility condition by leveraging the Helmholtz-Hodge decomposition. This creates a two-step time-stepping procedure for solving the Navier-Stokes equations, where an intermediate velocity field is first computed with the pressure omitted. Then, the intermediate velocity field is projected onto the divergence-free subspace to recover the incompressible velocity at the new time step. In this procedure, the pressure is never directly computed but can be reconstructed from the curl-free component of the intermediate velocity field (see \eg \ref{sec:pressure_reconstruction}). The decoupling of the momentum equation and the incompressibility condition via this projection method introduces temporal splitting errors that manifest on the boundary condition and can limit the order of accuracy of numerical methods in both space and time. However, the development of high-order projection methods has been an active area of research for the last three decades, and a number of standard approaches exist (see \cite{orszag1986boundary,karniadakis1991high,brown2001accurate,GUITTET2015215,THEILLARD201991,mancini2022projection}).

The overall stability of the projection method relies on the stability of each step. In the first step, the computation of the intermediate velocity field requires the solution of an advection-diffusion problem. The stability of this step is dictated by the properties of the temporal integration scheme and is thus easily assessed. The stability of the projection step, however, is linked to the spatial discretization of the gradient, divergence, and Laplacian operators. The natural way to ensure that the projection step is stable is to preserve the analytical properties of these operators and impose the divergence-free constraint at the discrete level. Unfortunately, this creates complications when using a collocated arrangement. For example, the discrete Laplacian is no longer a composition of the discrete gradient and divergence when using collocated variables with standard central difference formulas on Cartesian grids. These complications can be overcome by using interpolation procedures (\eg see \cite{minion1996projection}), but these procedures can themselves become challenging in the presence of irregular boundaries or adaptive grids. 

An alternate approach is to use a staggered layout, such as the Marker and Cell (MAC) representation introduced by Harlow and Welch \cite{harlow1965mac}, and illustrated in Figure \ref{fig:mac_gomez_us}. In this staggered arrangement, pressure data is located at the cell centers, and the velocities are located at the cell faces. In this context, the analytical properties of the operators are preserved, and the projection step enforces the divergence-free constraint with an orthogonal projection of the intermediate velocity field onto the divergence-free space. For this reason, the MAC layout has long been recognized as an ideal choice for solving the incompressible Navier-Stokes equations \cite{goda1979projection, kim1985application, van1986second}.

A different approach is to relax the requirement that the operators are preserved exactly at the discrete level and approximate the projection step. This notion of using an approximate projection was first introduced by Almgren \etal \cite{almgren1996approxproj} as a means of circumventing the numerical difficulties associated with exact discrete projections. In their approach, the discrete Laplacian operator is consistent, but not exactly the composition of the divergence and gradient, and the divergence of the resulting velocity field is only approximately zero up to the second order in the mesh spacing. This approximate projection uses velocities collocated at the cell centers with pressure data stored at the nodes. This design choice provides a symmetric discretization and generates a well-behaved linear system suitable for a multigrid solver. Additionally, the collocation of the velocities greatly simplifies the application of high-order upwind techniques and enables the overall time-stepping algorithm to achieve second-order convergence in both space and time. For a detailed analysis of approximate projection methods, we refer the reader to \cite{almgren2000approximate} and the references therein, which cover the use of approximate projections for a variety of incompressible flow problems and the combination of approximate projection methods with adaptive spatial and temporal meshes.

The inherent multiscale nature of the incompressible Navier-Stokes equations calls for adaptive mesh refinement (AMR) to optimize computation and memory overhead. Typically, only small localized regions of the domain need high grid resolution, for example, where high vorticity is present, such as in boundary layers or vortex cores. One of the earliest examples of using AMR is from Berger and Oliger \cite{berger1984adaptive}, where finer grids were adaptively placed over a coarse grid covering the domain. In \cite{almgren1997cartesian}, Almgren \etal combined an approximate projection method with an adaptive refinement strategy to solve the three-dimensional variable density incompressible Navier-Stokes equations. This block-structured approach was comprised of a nested hierarchy of logically-rectangular girds with refinement in both space and time. This AMR approach has been extended to numerous applications, and we refer the interested reader to the survey of \cite{dubey2014survey} and to \cite{schornbaum2018amr} for more details. For a modern block-structured AMR framework, we refer the reader to \cite{zhang2019amrex}.

Tree-based approaches to AMR \cite{Samet1988AnOO} are an alternative to block-structured AMR and the strategy used herein. They combine efficiency with simplicity by using recursive splitting schemes with non-overlapping regions. One of the earliest uses of graded octrees (a tree in which the size ratio between adjacent cells is at most two) for solving the incompressible Euler equations can be seen from Popinet in \cite{popinet2003gerris}, where finite volume discretizations are used on collocated cell-based quantities. This work uses the same approach as Almgren \etal in \cite{almgren1996approxproj} to approximate the Laplacian but also requires an approximation for the gradient operator. This leads to a nonuniform stencil and a non-symmetric system of linear equations for the pressure. In \cite{losasso2004simulating}, Losasso \etal proposed a symmetric discretization of the Poisson equation on octrees for free surface flows. This discretization resulted in a symmetric linear system, solved using a standard preconditioned conjugate gradient method. Their approach was only first-order accurate in the case of a non-graded mesh but was later extended to second-order accuracy \cite{losasso2006secondorder}, and later employed in the context of single phase \cite{Batty2007,GUITTET2015215,EGAN2021110084}, multiphase \cite{THEILLARD201991}, and free surface \cite{Enright2003,Ryoichi2020} applications. 

In \cite{min2006second}, Min \etal developed a collocated projection method for the incompressible Navier–Stokes equations on non-graded adaptive grids.  Their method utilized the Poisson solver of \cite{min2006supra}, which produced a non-symmetric but diagonally dominant linear system where the gradients of the solution were also second-order accurate. To ensure stability, instead of using an approximate projection, the authors manually enforced the orthogonality property between the divergence-free velocity field and the gradient of the Hodge variable when computing the updated velocity field. While this orthogonalization procedure results in an exact projection, this method is only stable if the normal velocity is null along the boundary. 

For a more versatile approach, Guittet \etal \cite{GUITTET2015215} developed a stable 
projection method for the incompressible Navier-Stokes equations on non-graded octree grids using a MAC layout. While this layout is a natural framework for a stable projection, the poor alignment between the pressure and velocity variables necessitated the use of complicated discretizations for the momentum equation and the projection step. The authors had to develop a Voronoi-based finite volume approach to treat the viscous terms implicitly and an expensive least squares interpolation was required to implement the semi-Lagrangian scheme for the advective terms. Though complicated discretizations were required, this framework was later extended to simulate active \cite{C6SM01955B,THEILLARD2019108841} and interfacial \cite{CLERETDELANGAVANT2017271,THEILLARD201991} flows.

An alternative staggered layout was proposed by Gomez \etal \cite{GOMEZ2019478}, in which the pressure was stored at the nodes and velocity stored along the edges (see Figure \ref{fig:mac_gomez_us}). This configuration of the staggered layout offers several advantages over the traditional MAC layout. The pressure gradient and velocity are naturally aligned, simplifying the construction of finite difference operators and granting higher accuracy. Unfortunately, this staggered layout, like the MAC configuration, requires multiple data structures and solvers. 

The method we present here is entirely collocated at the nodes of an arbitrary octree. It is carefully designed to achieve stability while reproducing solid interfaces and accompanying boundary layers with high fidelity. The overall stability of our method relies on the properties of our projection operator, which we analyze by relating the staggered and collocated approaches. We use a second-order semi-Lagrangian Backward Difference Formula scheme to update the momentum equation, treating the advective term explicitly and the viscous one implicitly. Arbitrary interface and boundary conditions are treated in a sharp manner using a level-set representation and the hybrid Finite-Volume/Finite-Difference discretizations from \cite{THEILLARD2019108841}. 

This manuscript is organized as follows. In Section \ref{sec:description}, we begin by recalling the incompressible Navier-Stokes equations and the general projection method used to solve them numerically.  In Section \ref{sec:projection}, we present the principle result of this work: our collocated projection operator on Cartesian grids. We prove its stability on uniform grids, construct a sufficient stability condition for adaptive grids, and provide numerical evidence of its stability for various boundary conditions and grid configurations. In Section \ref{sec:NSsolver}, we integrate our projection operator into a complete solver for the incompressible Navier-Stokes equations and observe overall second-order convergence. In Section \ref{sec:examples}, we illustrate the robustness and efficiency of our solver by using it to simulate several common validation problems of incompressible fluid flows in two and three spatial dimensions. We conclude in Section \ref{sec:ccl}.

\section{Mathematical model}\label{sec:description}
\subsection{Incompressible Navier-Stokes equations}
\begin{figure}[t!]
    \centering
    \includegraphics[width = .45\textwidth]{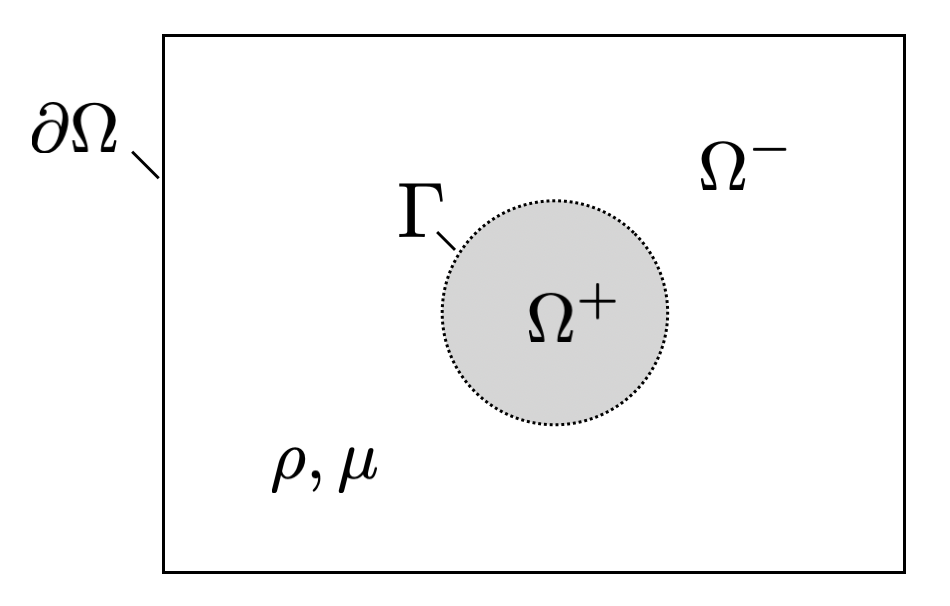}    
    \caption{Computational domain shown in two dimensions. The fluid domain, $\Omega^-$, is enclosed by the domain boundary, $\partial \Omega$, and the interface, $\Gamma$. Fluid properties, $\rho$ and $\mu$, are constant throughout the fluid domain. An arbitrary solid domain, $\Omega^+$ is shown as a shaded region.}
    \label{fig:domain}
\end{figure}

We consider a fluid set in a computational domain $\Omega = \Omega^+ \cup \Omega^- \subset \mathbb{R}^{2,3}$, where $\Omega^-$ represents the fluid phase and $\Omega^+$ represents all solid objects present in the fluid (see Figure \ref{fig:domain}), whose dynamics are modeled by the incompressible Navier-Stokes equations
\begin{align}
   \rho \left ( \frac{\partial \vectorr{u}}{\partial t} + \vectorr{u} \cdot \nabla \vectorr{u}\right ) &= - \nabla p +  \mu \Delta \vectorr{u} + \vectorr{f} \indomain, 
\label{eq:NS_mom} \\
    \nabla \cdot \vectorr{u}  &= 0  \indomain,
\label{eq:NS_mass}
\end{align}
where $\vectorr{u}$ is the fluid velocity, $p$ is the pressure, $\rho$ is the constant density, $\mu$ is the constant viscosity, and $\vectorr{f}$ are external forces such as the gravitational force. We denote the boundary of $\Omega^-$ by $\Gamma$ and the boundary of the computational domain $\Omega$ by $\partial \Omega$.

\subsection{General projection method}

The classical projection method is a fractional-step scheme for solving the incompressible Navier-Stokes equations. In the first step, often referred to as the viscosity or advection step, we advance the velocity field $\vectorr{u}^n$ at time step $t_n$ to an intermediate velocity field $\vectorr{u}^*$ by solving the momentum equation \eqref{eq:NS_mom}  where the pressure term is omitted. For example, using a standard first-order semi-explicit scheme, $\vectorr{u}^*$ is constructed as the solution of
\begin{equation}
    \rho \left (\frac{\vectorr{u}^*-\vectorr{u}^n}{\Delta t} + \vectorr{u}^n \cdot \nabla \vectorr{u}^n \right ) = \mu \Delta \vectorr{u}^* + \vectorr{f}. \label{eqn:discretizedmomentum}
\end{equation}
Next, we project $\vectorr{u}^*$ onto the divergence-free space to enforce the incompressibility condition \eqref{eq:NS_mass}. To do so, we use the Helmholz-Hodge decomposition to separate the intermediate velocity into curl-free and divergence-free components as follows:
\begin{equation}
\label{eq:hodgedecomp}
    \vectorr{u}^* = \vectorr{u}^{n+1} + \nabla \phi,
\end{equation}
where  $\vectorr{u}^{n+1}$ is the divergence-free velocity
field at time $t_{n+1}$ and $\phi$ is the Hodge variable. Taking the divergence of the above equation  \eqref{eq:hodgedecomp},  we obtain a Poisson equation for the Hodge variable,
\begin{equation}
    \Delta \phi = \nabla \cdot \vectorr{u}^*.
\end{equation}
Using the appropriate boundary conditions, the above equation is solved to compute the Hodge variable, and the divergence-free velocity is recovered as
\begin{equation}
    \vectorr{u}^{n+1} = \vectorr{u}^* - \nabla \phi.
    \label{eq:unp1projection}
\end{equation}
Note that equation \eqref{eq:unp1projection} can be rewritten as an operator applied to $\vectorr{u}^\ast$,
\begin{equation}
\vectorr{u}^{n+1} =\left(\operator{I} - \nabla \Delta ^{-1} \nabla \cdot\right) \vectorr{u}^*. 
\end{equation}
Thus, we define the generic projection operator ${P}$ as, 
\begin{equation}
{P} = \operator{I} - \nabla \Delta ^{-1} \nabla \cdot. 
\label{eq:pdiffopdef}
\end{equation}
\subsection{Properties of the projection operator} 
\label{sec:generalstuff}
The analytic operator ${P}$, defined by equation  \eqref{eq:pdiffopdef}, is indeed a projection (\ie ${P}^2={P}$) as
\begin{align}
{P}^2=\left(\operator{I} - \nabla \Delta ^{-1} \nabla \cdot\right)^2 &= \operator{I} - 2\nabla \Delta ^{-1} \nabla \cdot
+\nabla \Delta ^{-1} \nabla \cdot \nabla \left(\nabla \cdot \nabla\right)^{-1} \nabla \cdot,\\
&= \operator{I} - \nabla \Delta ^{-1} \nabla \cdot.
\end{align}
The projection property ensures that the projected field is exactly divergence-free and, thus, that the incompressibility condition is exactly satisfied. Moreover, since the gradient and divergence are the negative transpose of each other ($\nabla^T = -\nabla \cdot$), the projection is symmetric $P^T =P$, orthogonal \ie
\begin{align}
    \norm{\vectorr{u}^{n+1}}^2 =    \norm{P\vectorr{u}^*}^2 &= \norm{\vectorr{u}^*}^2 -2<\vectorr{u}^*| \nabla\Phi> + \norm{\nabla\Phi}^2, \\
   &= \norm{\vectorr{u}^*}^2 +2<\nabla \cdot\vectorr{u}^*| \Phi> + \norm{\nabla\Phi}^2,\\
   &= \norm{\vectorr{u}^*}^2 +2<\nabla \cdot\nabla\Phi| \Phi> + \norm{\nabla\Phi}^2,\\
    &= \norm{\vectorr{u}^*}^2 - \norm{\nabla\Phi}^2, 
\end{align}
and therefore contracting,
\begin{equation}
\norm{P\vectorr{u}^*} \leq \norm{\vectorr{u}^{*}}.
\end{equation}
In the discrete case, if the composition and negative transpose properties are preserved, $P$ remains contracting, and, provided that the temporal integrator is stable, the entire update from $t_n$ to $t_{n+1}$ is stable.

\subsection{Boundary conditions}
\label{sec:bcs}
We focus on single-phase flows around solid objects; hence, the possible boundary conditions are no-slip at the walls and interface and possible inflow/outflow boundary conditions at the walls of the computational domain. The no-slip boundary condition is a Dirichlet boundary condition on the velocity, $\vectorr{u}|_{\Gamma},$  where $\vectorr{u}|_{\Gamma}$ is zero if the interface is static and equal to the velocity of the interface if it is moving. This, along with equation (\ref{eq:hodgedecomp}), implies that the boundary condition of the intermediate velocity field $\textbf{u}^*$ is
\begin{equation}
\vectorr{u}^*|_{\Gamma} = \vectorr{u}|_\Gamma + \nabla\phi^{n+1}|_\Gamma .
\label{eq:bcustar}
\end{equation}
When solving for the Hodge variable in the projection step, we prescribe the homogeneous Neumann boundary conditions, \ie
\begin{equation}
(\nabla\phi \cdot \vectorr{n})|_{\Gamma} = 0,
\end{equation}
where $\vectorr{n}$ is the outward-facing normal vector to the boundary. Note that the value of $\phi^{}$ is not computed when finding the intermediate velocity field $\textbf{u}^*.$ To keep the viscosity and projection decoupled, we can use the Hodge variable at the previous time step, $\phi^n,$ to serve as an initial guess to the correct boundary condition of $\vectorr{u}^*$ and iterate until we have convergence in both $\vectorr{u}^*$ and $\phi$. In \cite{GUITTET2015215,EGAN2021110084,mancini2022projection}, it was observed that using this initial guess led to very few iterations being needed to reach convergence. Other iterative corrections can be designed, (see for example \cite{THEILLARD2019108841}). The $\nabla \phi^{n+1}$ term in the above \eqref{eq:bcustar} should be seen as a splitting error. It can be shown to be first-order in time. Even in the absence of iterative correction, the presence of this term will only introduce converging errors on the boundary conditions \cite{brown2001accurate}. 

Potential flux boundary conditions on the walls of the computational domain $\partial \Omega$ correspond to Neumann boundary conditions. This means that the boundary condition for the intermediate velocity field $\vectorr{u}^*$ is
\begin{equation}
    (\nabla\vectorr{u}^*\cdot\vectorr{n})|_{\partial \Omega} = (\nabla\vectorr{u}\cdot \vectorr{n})|_{\partial\Omega} + (\nabla\nabla\Phi \cdot \vectorr{n})_{\partial \Omega}.
\end{equation}
The boundary condition for the Hodge variable is a Dirichlet boundary condition. Again the last term in the above equation should be seen as a splitting and converging ($\mathcal{O}(\Delta t)$) error. In this case, an iterative correction for this term would require differentiating the Hodge variable twice, which may produce noisy and inaccurate results. Therefore, this term is often disregarded. If boundary conditions are prescribed on the pressure directly, the corresponding boundary conditions on the Hodge variable are obtained from the pressure reconstruction formula.

\subsection{Pressure reconstruction}
\label{sec:pressure_reconstruction}
In our formulation of the projection method, pressure is never directly computed. If desired, we can compute it using 
\eqref{eqn:discretizedmomentum} and 
\eqref{eq:hodgedecomp} by asking that the discrete momentum equation must be satisfied for $\vectorr{u}^{n+1}$ and $p$, and obtaining
\begin{equation}
    \rho \left (\frac{\vectorr{u}^{n+1} - \vectorr{u}^n}{\Delta t} + \vectorr{u}^n \cdot \nabla \vectorr{u}^n \right ) = - \nabla p + \mu \Delta \vectorr{u}^{n+1}.
\end{equation}
Since $\vectorr{u}^{n+1} = \vectorr{u}^* - \nabla \phi$ from the Helmholtz-Hodge decomposition
\begin{equation}
    \rho \left (\frac{\vectorr{u}^* - \nabla \phi - \vectorr{u}^n}{\Delta t} + \vectorr{u}^n \cdot \nabla \vectorr{u}^n \right ) = - \nabla p + \mu \Delta (\vectorr{u}^* - \nabla \phi),
\end{equation}
and because $\vectorr{u}^*$ is defined as the solution of \eqref{eqn:discretizedmomentum}, this simplifies into
\begin{equation}
-\rho  \frac{\nabla \phi}{\Delta t} = -\nabla p - \mu \Delta \nabla \phi \text{,}
\end{equation}
or equivalently, 
\begin{equation}
\nabla p = \frac{\rho}{\Delta t}\nabla  \phi - \mu \nabla\Delta \phi ,  
\end{equation} 
and so, up to an additive constant, we arrive at our pressure reconstruction formula
\begin{equation}
p = \frac{\rho}{\Delta t} \phi - \mu \Delta \phi.
\end{equation}
We point out that using a different time integrator may alter the above formula (see \eg \cite{GUITTET2015215,THEILLARD201991}).

\section{Node-based projection operator on Cartesian grids } \label{sec:projection}

In this section, we introduce our collocated projection method, analyze its properties on periodic uniform and adaptive Cartesian grids, and verify the stability of the method numerically. 
We find that the composition projection property is lost due to the collocation, and thus our operator is not a projection. However, we see that our operator can be iterated to recover the canonical projection onto the divergence-free space. We derive a technical sufficient condition for these iterations to converge and, in the periodic uniform case, are able to prove that our operator is, in fact, contracting. 

The convergence analysis relies on the connection between the staggered and collocated frameworks. To this effect, we show that the collocated projection can be related to the staggered one through interpolation procedures and leverage important facts about the discrete projection operator on staggered grids. The interpolations are solely introduced for the purpose of proving the stability of the collocated projection; they do not appear in the numerical implementation. 

\subsection{Definitions and key observation} \label{sec:pres_and_gen_prop}
 Our nodal projection operator $\operator{P}_N$ is defined from the nodal gradient $\operator{G}_N$, divergence $\operator{D}_N$ and Laplacian $\operator{L}_N$ as
\begin{equation}
    \operator{P}_N = \operator{I} - \operator{G}_N \operator{L}_N^{-1} \operator{D}_N,
    \label{eq:projoperdef}
\end{equation}
where $\operator{I}$ is the appropriate identity operator and $\operator{L}_N^{-1}$ is the operator that for any vector $X$, returns $Y$ the solution  of $\operator{L}_NY =X$ with problem-dependent homogeneous boundary conditions \footnote{If the boundary conditions are chosen so that the solution is not uniquely defined, we will enforce that the solution is zero at one specific location.}. 
\begin{figure}[b!]
\centering
    \includegraphics[width = .95\textwidth]{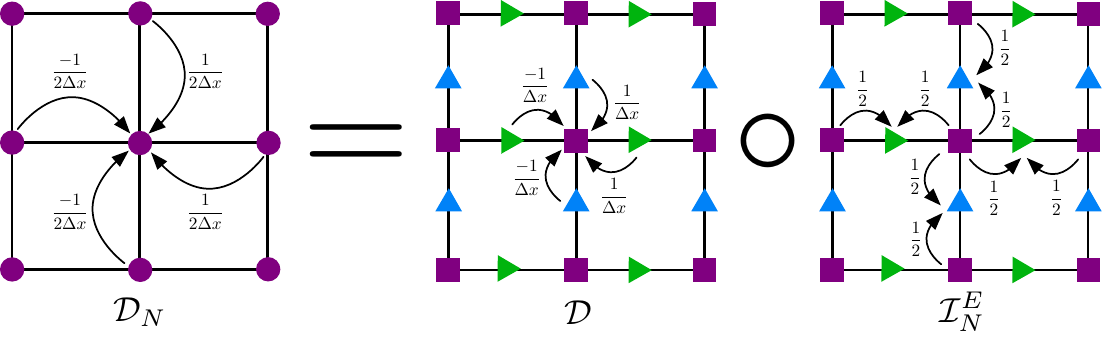}    
    \caption{Connection between the discrete divergence on collocated ($\operator{D}_N$) and staggered grids ($\operator{D}$). The arrows represent the finite difference coefficients used in each contributing term.}
    \label{fig:Dn=DI}
\end{figure}
Again, the nodal Laplacian is not the composition of the nodal divergence and nodal gradient (\ie $\operator{L}_N \ne \operator{G}_N \operator{D}_N$) and our nodal projection operator is not a true projection (\ie $\operator{P}_N^2\neq \operator{P}_N$). This means that the projected velocity field is not guaranteed to be divergence-free (\ie $\operator{D}_N\operator{P}_NX\neq 0$). However, we can show that our projection operator only preserves the incompressible modes of the velocity field. Therefore, if we iterate our projection operator and if the iterated operator converges, then it must converge to the canonical orthogonal projection. 

To demonstrate this, we assume that $\operator{P}^k_N\rightarrow\operator{P}^\infty_N$ as $k\rightarrow \infty$, and remark that if $\lambda_i$ are the eigenvalues of $\operator{P}_N$, the eigenvalues of $\operator{P}^\infty_N$ are  $\lim_{k\rightarrow \infty} \lambda_i^k$, and since they are finite, they can only be 1 or 0. In addition, the eigenvectors of $\operator{P}^\infty_N$ corresponding to the eigenvalue 1 are the eigenvectors of $\operator{P}_N$ for the same eigenvalue. Therefore $\operator{P}^\infty_N$ is the projection on the eigenspace of $\operator{P}_N$ corresponding to the eigenvalue 1. Clearly, any divergence-free modes $X$ belong to this eigenspace since
\begin{equation}
    \operator{P}^\infty_NX=\lim_{k\rightarrow\infty}\operator{P}^k_NX = X.
\end{equation}
Now, if $X$ is an eigenvector associated to the eigenvalue $1$, $\operator{P}_N {X} = {X}$, and so
\begin{equation}
 \operator{G}_N \operator{L}_N^{-1} \operator{D}_N X =0,
\end{equation}
meaning that $\operator{L}_N^{-1} \operator{D}_N {X}$ must be constant. Assuming that a homogeneous Dirichlet boundary condition is enforced somewhere, this constant must be zero 
\begin{equation} 
\operator{L}_N^{-1} \operator{D}_N X = 0 ,
\end{equation}
and thus $ \operator{D}_N X = 0$, indicating that $X$ is incompressible. The eigenspace for the eigenvalue 1 is the incompressible space, and therefore $P^\infty_N$, if it exists, is the canonical projection on the incompressible space. 

In the subsequent sections, we construct a sufficient condition for this limit to exist. However, we note that in practice, this limit may not be reached as only a small number of iterations will be performed, and the advanced velocity field may not be exactly divergence-free. These deviations can be kept arbitrarily small by adjusting the stopping criteria, and in practice, a few iterations are needed to achieve a reasonably small error tolerance. In fact, a single iteration may be enough to obtain an overall second-order solution (see Section \ref{sec:det_kmax}).

 \subsection{Uniform Cartesian grids}

\subsubsection{Periodic domains}

\begin{figure}[t!]
    \centering
     \includegraphics[width = .95\textwidth]{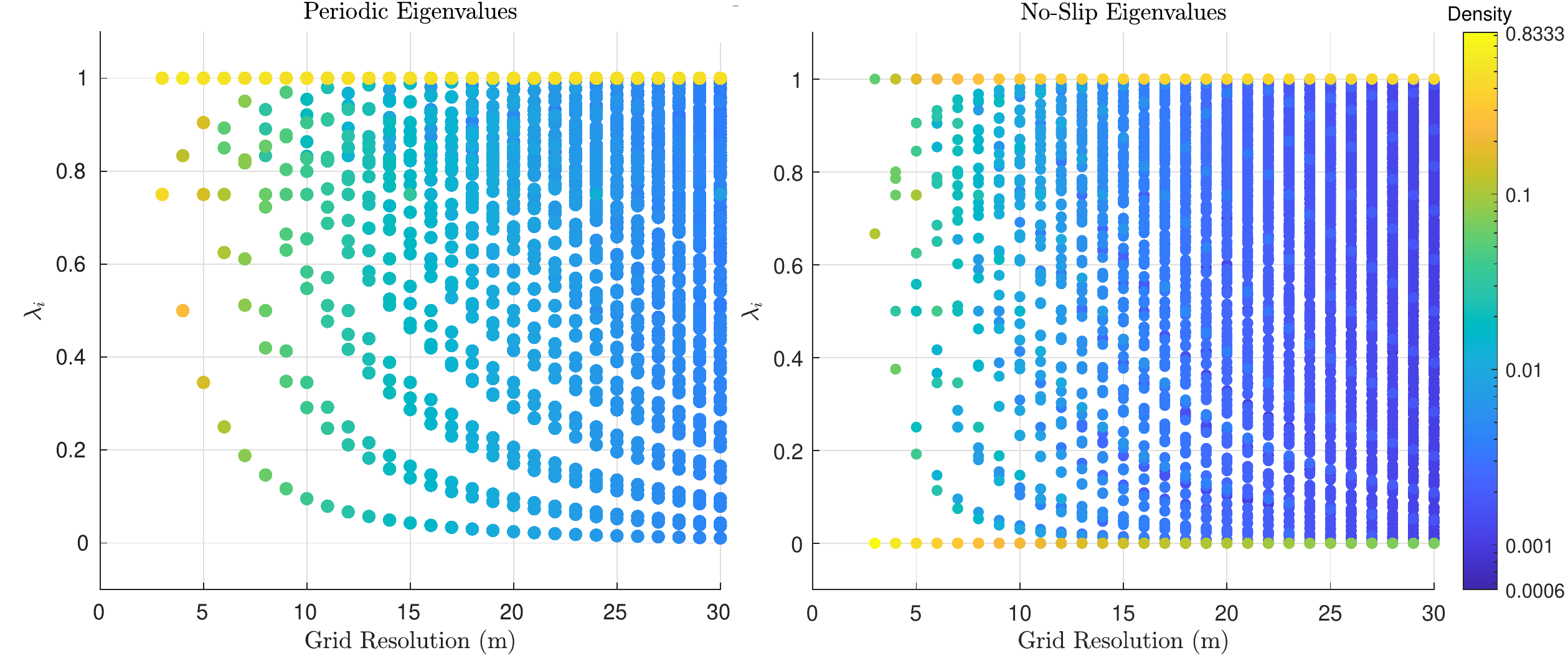}
    \caption{Spectrum of the projection operator with periodic boundary conditions (left) and no-slip boundary conditions (right) on a uniform grid. We observe that the eigenvalues are in the interval [0,1], which suggests that the operator is stable. The density of the eigenvalues for each grid resolution is shown. }
    \label{fig:projectioneigs_density}
\end{figure}

In this case, we define  $\operator{G}_N$, $\operator{L}_N$, $\operator{D}_N$ to be the standard periodic second-order finite difference gradient, Laplacian, and divergence. As Figure \ref{fig:Dn=DI} illustrates, the nodal divergence operator $\operator{D}_N$ can be seen as the composition of the staggered discrete divergence $\operator{D}$ with the linear interpolation operator from the set of nodes $N$ to the set of edges $E$ $\operator{I}_N^E:N\rightarrow E$
\begin{equation}
\operator{D}_N = \operator{D} \operator{I}_N^E.
\end{equation}
Similarly, the discrete nodal gradient $\operator{G}_N$ is related to the staggered operator and the interpolation from the edges to the nodes $\operator{I}_E^N$  as
\begin{equation}
\operator{G}_N = \operator{I}_E^N \operator{G}.
\end{equation}
The staggered operators satisfy the negative transpose property
\begin{equation}
    \operator{D}^T =-\operator{G},
\end{equation}
and the two interpolations are transpose of each other
\begin{equation}
    (\operator{I}^N_E)^T =\operator{I}_N^E.
\end{equation}
Finally, we observe that the nodal Laplacian $\operator{L}_N$ is identical to the staggered one $\operator{L}$
\begin{equation}
\operator{L}_N =\operator{L}.
\end{equation}
Using these observations, we can rewrite the projection operator \eqref{eq:projoperdef} as
\begin{equation}\label{eq:projoperator_interp}
    \operator{P}_N = \operator{I} - I_E^N\operator{G}\operator{L}^{-1}\operator{D}I_N^E.
\end{equation}

All these ingredients allow us to follow the contraction's proof  (see Section \ref{sec:generalstuff}): the norm of the projected velocity is
\begin{equation}
    \norm{\operator{P}_NX}^2 = \norm{X}^2 -2<X|I_E^N\operator{G}\operator{L}^{-1}\operator{D}I_N^EX> + \norm{I_E^N\operator{G}\operator{L}^{-1}\operator{D}I_N^EX}^2,
\end{equation}
setting $\operator{L}^{-1}\operator{D}I_N^EX=\Phi$, and 
using the transpose properties of the divergence and interpolation we get
\begin{equation}
    \norm{\operator{P}_NX}^2 = \norm{X}^2 +2<\operator{D}I_N^EX|\Phi> + \norm{I_E^N\operator{G}\Phi}^2,
\end{equation}
which is identical to 
\begin{equation}
    \norm{\operator{P}_NX}^2 = \norm{X}^2 +2<\operator{L}\operator{L}^{-1}\operator{D}I_N^EX|\Phi> + \norm{I_E^N\operator{G}\Phi}^2,
\end{equation}
giving us
\begin{equation}
    \norm{\operator{P}_NX}^2 = \norm{X}^2 +2<\operator{L}\Phi|\Phi> + \norm{I_E^N\operator{G}\Phi}^2,
\end{equation}
and since $\operator{L}=\operator{D}\operator{G}$ and $\operator{G}=-\operator{D}^T$
\begin{equation}
    \norm{\operator{P}_NX}^2 = \norm{X}^2 -2<\operator{G}\Phi|\operator{G}\Phi> + \norm{I_E^N\operator{G}\Phi}^2,
\end{equation}
we have
\begin{equation}
    \norm{\operator{P}_NX}^2 = \norm{X}^2 -2\norm{\operator{G}\Phi}^2+ \norm{I_E^N\operator{G}\Phi}^2,
\end{equation}
and so
\begin{equation}
    \norm{\operator{P}_{N}X}^2 \leq \norm{X}^2 -\left(2- \norm{I_E^N}^2\right)\norm{\operator{G}\Phi}^2.
\end{equation}
Thus, if $\norm{I_E^N}<\sqrt{2}$, we have $ \norm{\operator{P}_{N}X}\leq \norm{X}$ and $\operator{P}_{N}$ is contracting, and therefore converging. Here, the interpolations are averaging operators, so their $L^\infty$-norms are less than or equal to one, and the projection is contracting. 

In Figure \ref{fig:projectioneigs_density}, we display the spectrum of the projection operator, computed directly from the discrete numerical operator and for increasing grid resolution. Since the operator is contracting, we expect its spectrum to lie in the interval $[0,1]$, and indeed it does. We color the eigenvalues by the density of the eigenvalue for that grid resolution. Note that for all grid resolutions, the eigenvalue with the highest density is the $\lambda = 1$ eigenvalue, which corresponds to an eigenvector with an incompressible mode.

\subsubsection{No-slip boundary conditions}
For a more realistic  context, we turn our attention to no-slip boundary conditions. On rectangular domains, this means that we set the $x$-velocity $u = 0$ along the horizontal walls, the $y$-velocity $v = 0$ along the vertical walls, and the normal gradient of the Hodge variable to be zero on all walls. Currently, we do not have proof of the stability of this operator, as a general formulation of this operator with these boundary conditions, unlike in the periodic case, is not easily known. Instead, we reconstruct the numerical operator by projecting the canonical basis vectors of $\mathbb{R}^n$, where $n$ is the dimension of the square matrix operator, with the boundary condition set beforehand. The resulting spectrum for increasing grid resolution is depicted in Figure \ref{fig:projectioneigs_density}, and is contained in $[0,1]$. The addition of the boundary condition introduces an extra constraint on the velocity field, causing the projection to be even more restrictive and thus the admissible incompressible modes (\ie eigenvector associated to $\lambda=1$) to be rarer. In addition, the no-slip boundary condition introduces a 0 eigenvalue to the operator, which corresponds to the canonical basis vectors that are non-zero only along the walls with the prescribed boundary condition.

\subsection{Extension to Adaptive Grids}
The negative transpose property is lost on adaptive grids, so the projected velocity norm cannot easily be bounded. Instead, we focus on the spectrum of the projection, and since we know that the incompressible space is the eigenspace associated with the eigenvalue 1 (see Section \ref{sec:pres_and_gen_prop}), we only need to find a sufficient condition for all eigenvalues associated to compressible eigenvectors to be strictly less than 1 in absolute value. In other words, a condition under which our projection is critically stable.

We start our study by formally defining all involved operators and proving the stability of the staggered projection. We then leverage this property to find a sufficient convergence condition. For the sake of clarity, this entire theoretical study is done on two-dimensional periodic grids. The extension to three dimensions is tedious but straightforward. 

\begin{figure}[h]
    \centering   
    \includegraphics[width = .85\textwidth]{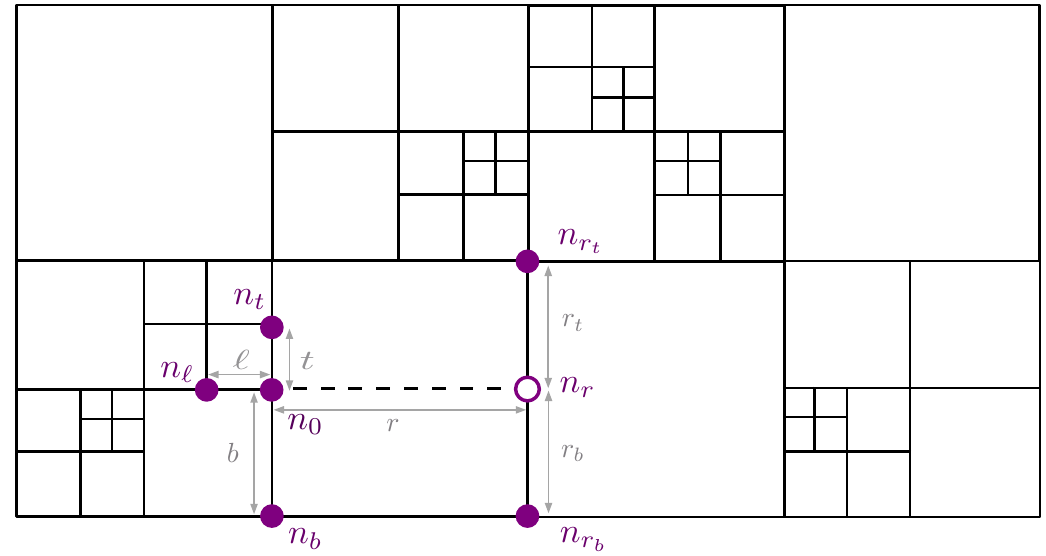}
    \caption{Finite difference discretization on  quadtree grids. Here, node $n_0$ has no direct neighbor to the right, and thus a ghost node $n_r$ (\includegraphics[height=0.015\textwidth]{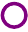}) must be constructed using the existing neighboring nodes (\includegraphics[height=0.015\textwidth]{figures/purpledisk.pdf}). Standard central discretizations can then be constructed using this ghosted neighborhood.}
    \label{fig:laplaciannodes}
\end{figure}

\begin{figure}[h]
    \centering
    \includegraphics[width = .95\textwidth]{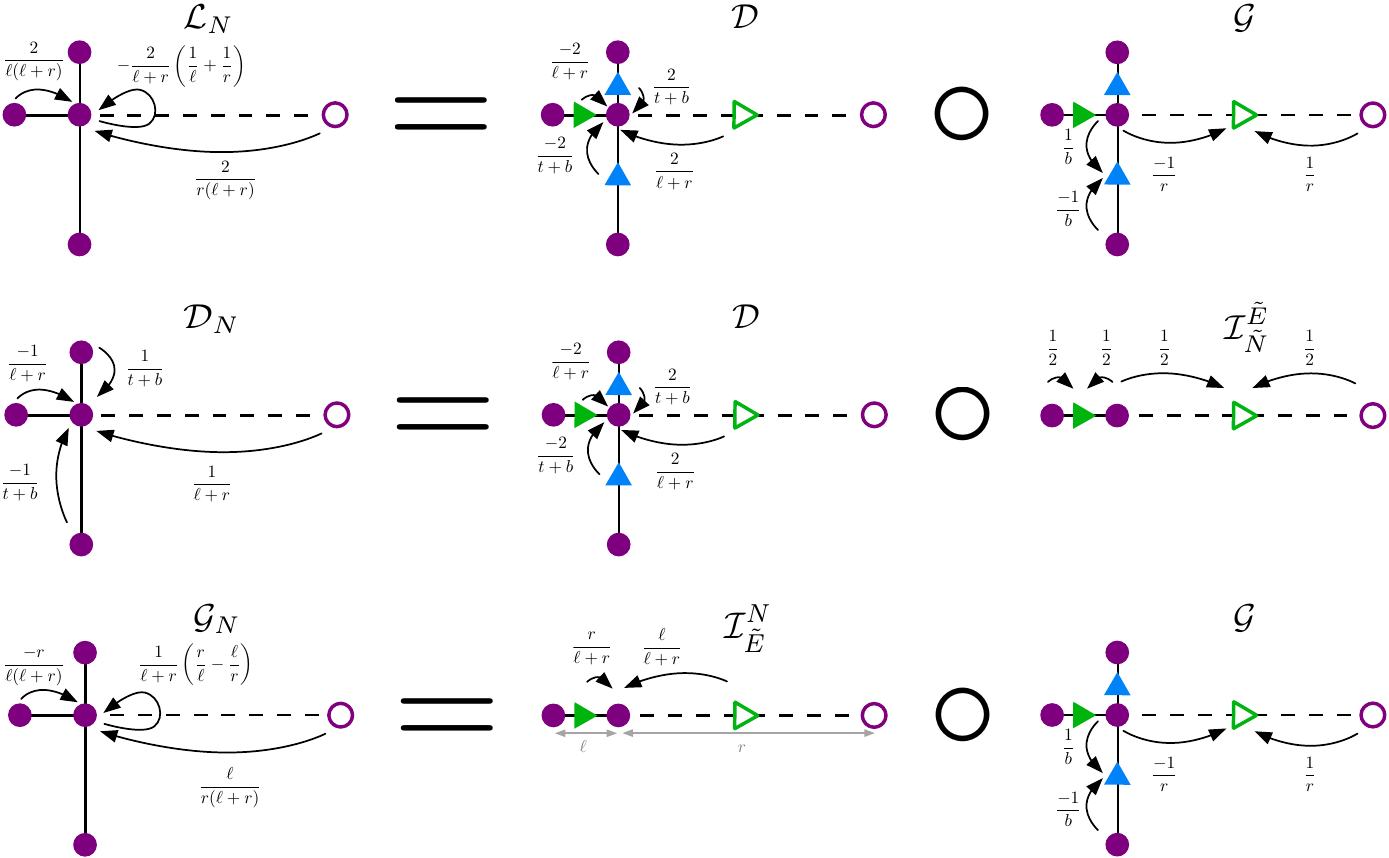}
    \caption{Connections between the staggered and collocated operators on periodic quadtree grids. The nodal Laplacian $\operator{L}_N$ is the composition of the staggered divergence $\operator{D}$ and gradient $\operator{G}$. The nodal divergence $\operator{D}_N$ and gradient $\operator{G}_N$ are related to their staggered counterpart through linear interpolation. Ghost variables are depicted as empty symbols (\includegraphics[height=0.015\textwidth]{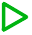} and \includegraphics[height=0.015\textwidth]{figures/emptycircle.pdf}). As in Figure \ref{fig:laplaciannodes}, arrows represent the finite difference coefficients. For the gradient and Laplacian, only the x-contribution is illustrated. The y-contribution is computed similarly.}
    \label{fig:operators_quadtree}
\end{figure}

\subsubsection{Construction of the discrete differential operators}

To discretize our differential operators, we follow the work of Min and Gibou \cite{min2006supra}. The central idea is to define ghost values at T-junctions (often called hanging nodes) to circumvent the lack of direct neighbors. This is done using standard Taylor analysis and derivative approximation. For example, referring to Figure \ref{fig:laplaciannodes}, node $n_0$ does not have a direct neighbor to the right, so we introduce a ghost node $n_r$ on the face delimited by nodes $n_{rt}$ and $n_{rb}$. For any nodal quantity $\phi$, sampled at the existing nodes, we can calculate a third-order accurate ghost value $\phi_r$ using the information at $n_0$, at its direct neighbors in all three other directions (\ie $n_l,n_t,n_b$), and at the neighboring nodes $n_{rt}$ and $n_{rb}$ as 
\begin{equation}
\phi_r = \frac{r_b\phi_{r_t} +  r_t\phi_{t_b}}{r_t+r_b} - \frac{r_t r_b}{t+b }\left( \frac{\phi_t-\phi_0}{t}- \frac{\phi_0-\phi_b}{b}\right).
\label{eq:ghostnode_int}
\end{equation}
Now, we should see each node as having four direct neighbors, one in each direction, with at most one ghost neighbor. Using these new neighborhoods, we can construct discretizations for the standard differential operators, such as the nodal Laplacian $\operator{L}_N$, divergence $\operator{D}_N$ and gradient $\operator{G}_N$. As Figure \ref{fig:operators_quadtree} illustrates,
these operators evaluated at any node $n_0$ for any given sampled Hodge variable $\phi$ and any given velocity field $\vectorr{u}=\left({u},{v}\right)$ are expressed in terms of their neighborhood information as
\begin{align}
\operator{L}_N \phi\big|_0 &= \frac{2}{r+\ell}\left(\frac{\phi_r-\phi_0}{r} - \frac{\phi_0 - \phi_\ell}{\ell}\right) +\frac{2}{t+b}\left(\frac{\phi_t-\phi_0}{t} - \frac{\phi_0 - \phi_b}{b}\right), \\
\operator{D}_N(u,v)\big|_0 &= \frac{u_r-u_\ell}{r+\ell} + \frac{v_t-v_b}{t+b},\\
\operator{G}_N\phi\big|_0 &= \left(
\frac{\ell }{r+\ell} \frac{\phi_r-\phi_0}{r}+ \frac{r}{r+\ell}\frac{\phi_0-\phi_\ell}{\ell}, \quad \frac{b}{t+b}\frac{\phi_t-\phi_0}{t} + \frac{t}{t+b}\frac{\phi_0-\phi_b}{b} \right).
\end{align}

Note that our definitions of the Laplacian and gradient are identical to those in previous studies (see \eg \cite{min2006supra, min2006second, THEILLARD2013430,  Theillard2013AMM, THEILLARD2019108841, KUCHEROVA2021110591}), and in particular, they have been shown to yield supra convergence, in the sense that both the solution and its gradient converge with second-order accuracy \cite{min2006supra}. The divergence operator, however, does not follow the construction proposed by Min and Gibou \cite{min2006supra,min2006second}. Here, we use the standard central difference formula to construct a first-order accurate approximation of the divergence operator instead of creating weighted averages of the forward and backward derivative and constructing a second-order centered difference approximation. It was observed in \cite{min2006second} that using the second-order divergence in a purely nodal context leads to an unstable projection operator.

To establish the connection between the staggered and collocated projection, we must formally define all discrete differential operators and the interpolations needed to bridge between them, as well as the spaces they act on. Using ghost nodes allows us to see these operators as natural extensions of their uniform counterparts.

\subsubsection{Interpolations: notations and definitions}
From the set of nodes and edges $N$ we define  $\tilde{N}$ as the set of all grid nodes and ghost neighbors (depicted as empty circles in Figure \ref{fig:operators_quadtree}), and refer to it as the set of ghosted nodes. For each ghost node, we define a ghost edge  (depicted as empty triangles in Figure \ref{fig:operators_quadtree}), centered between the ghost node and the corresponding hanging node. We denote by $\tilde{E}$ the set of ghosted edges, \ie grid, and ghost edges. As we did before, for any two discrete spaces $A$ and $B$, defined from the grid structure, we will define $\operator{I}_A^B:A\rightarrow B$ as the interpolation from $A$ to $B$.

As Figure \ref{fig:operators_quadtree} illustrates, we define the interpolation from the set of ghosted nodes to ghosted edges, $\operator{I}_{\tilde{N}}^{\tilde{E}} : \tilde{N} \to \tilde{E}$, and the interpolation from the set of ghosted edges to the set of nodes, $\operator{I}_{\tilde{E}}^{{N}} : \tilde{E} \to {N}$, as the standard linear interpolations from the nearest neighbors. We also define $\operator{I}_{N}^{\tilde{N}}  : N \to \tilde{N}$ as the ghost operator, which returns the value at the grid nodes from the interpolated value at the ghost nodes. Similarly, we define $\operator{I}_{E}^{\tilde{E}} : E \to \tilde{E}$.

We define the operator $\operator{I}_{\tilde{E}}^{E}$ as the canonical restriction from $\tilde{E}$ to $E$. The staggered gradient operator, $\operator{G} : \tilde{N} \to \tilde{E}$, is constructed using the standard second-order central difference scheme and the staggered divergence operator, $\operator{D} : \tilde{E} \to N$ is constructed using a weighted first-order central difference scheme (see Figure \ref{fig:operators_quadtree}). 

\subsubsection{Staggered projection on quadtree grids}

We start by observing that the nodal Laplacian is obtained by taking the divergence of the gradient computed for the ghosted values, \ie $\operator{L}=\operator{D}\operator{G}\operator{I}_{{N}}^{\tilde{N}}$, and so the projection operator  $\operator{P}:{E}\rightarrow{E}$ has the form
\begin{equation}
\operator{P} =\operator{I} - \operator{I}_{\tilde{E}}^{{E}}\operator{G} \operator{I}_{{N}}^{\tilde{N}}\left( \operator{D}\operator{G}\operator{I}_{{N}}^{\tilde{N}}\right)^{-1}\operator{D}\ghostinterpedges.
\end{equation}
Because $\operator{I}_{\tilde{E}}^{{E}}\ghostinterpedges=\operator{I}$ (i.e un-ghosting the ghosted edges' values returns the edges' values), we can rewrite the above definition by defining the operator $\operator{P}_{\tilde{E}}$ as 
\begin{equation}
\operator{P} =\operator{I}_{\tilde{E}}^{{E}}\underbrace{\left(\operator{I} - \operator{G} \operator{I}_{{N}}^{\tilde{N}}\left( \operator{D}\operator{G}\operator{I}_{{N}}^{\tilde{N}}\right)^{-1}\operator{D}\right)}_{\operator{P}_{\tilde{E}}}\ghostinterpedges.
\label{eq:pstagghost}\end{equation}
$\operator{P}_{\tilde{E}}$ should be interpreted as the projection on the space of ghosted edges. Squaring it, we see that 
\begin{align}
\operator{P}_{\tilde{E}}^2 &=\operator{I} - 2\operator{G} \operator{I}_{{N}}^{\tilde{N}}\left( \operator{D}\operator{G}\operator{I}_{{N}}^{\tilde{N}}\right)^{-1}\operator{D} + \operator{G} \operator{I}_{{N}}^{\tilde{N}}\left( \operator{D}\operator{G}\operator{I}_{{N}}^{\tilde{N}}\right)^{-1}\operator{D}\operator{G} \operator{I}_{{N}}^{\tilde{N}}\left( \operator{D}\operator{G}\ghostinterp\right)^{-1}\operator{D},\\
\operator{P}_{\tilde{E}}^2 &=\operator{I} - \operator{G} \operator{I}_{{N}}^{\tilde{N}}\left( \operator{D}\operator{G}\ghostinterp\right)^{-1}\operator{D},\\ 
\operator{P}_{\tilde{E}}^2 &=\operator{P}_{\tilde{E}},
\end{align}
and so $\operator{P}_{\tilde{E}}$ is indeed a projection\footnote{Note that, as it is the case in the continuous world, this propriety only relies on the fact that the staggered Laplacian is the composition of the divergence  and  gradient operators. We, for example, did not have to assume that these two operators are the negative transpose of each other.}.

\subsubsection{Stability of the collocated projection}
To connect the staggered and collocated projections, we start by recognizing that the nodal Laplacian $\operator{L}_N$ is the composition of the staggered divergence with the staggered gradient (see Figure \ref{fig:operators_quadtree})
\begin{equation}
\operator{L}_N = \operator{D}\operator{G}\ghostinterp.
\end{equation}
As it was the case on uniform grids, the nodal divergence and gradient remain related to their staggered counterparts through linear interpolations, only now these relationships involve the ghost interpolation operator $\ghostinterp$
\begin{align}
\operator{D}_N = \operator{D}\operator{I}_{\tilde{N}}^{\tilde{E}}\ghostinterp, \\
\operator{G}_N =\operator{I}_{\tilde{E}}^{{N}}\operator{G}\ghostinterp. 
\end{align}
Using the above expression, the nodal projection operator 
\begin{equation}
\operator{P}_N = \operator{I} - \operator{G}_N\left( \operator{L}_N\right)^{-1}\operator{D}_N \text{,}
\end{equation}
can be expressed as
\begin{equation}
\operator{P}_N = \operator{I} -\operator{I}_{\tilde{E}}^{{N}}\operator{G} \operator{I}_{{N}}^{\tilde{N}}\left( \operator{D}\operator{G}\operator{I}_{{N}}^{{\tilde{N}}}\right)^{-1}\operator{D}\operator{I}_{\tilde{N}}^{\tilde{E}}\ghostinterp.
\end{equation}
We consider a compressible eigenpair $\left(\lambda, X\right)$ such that 
\begin{equation}
\operator{P}_{N}X = \lambda X,  \qquad \operator{D}_N{X}\neq 0, \qquad \textrm{and} \qquad \lambda\neq 1. 
\end{equation}
We multiply the whole equation on the left by the operator $\operator{I}_{\tilde{N}}^{\tilde{E}}\ghostinterp$, and define   $Y =\operator{I}_{\tilde{N}}^{\tilde{E}}\ghostinterp X $, to obtain
\begin{equation}
\operator{I}_{\tilde{N}}^{\tilde{E}}\ghostinterp \operator{P}_{N}X = Y - \underbrace{\operator{I}_{\tilde{N}}^{\tilde{E}}\ghostinterp \operator{I}_{\tilde{E}}^{{N}}}_{\operator{I}_p}
\underbrace{\operator{G} \operator{I}_{{N}}^{\tilde{N}}\left( \operator{D}\operator{G}\operator{I}_{{N}}^{{N}}\right)^{-1}\operator{D}}_{\operator{Q}_{\tilde{E}}=\operator{I}-\operator{P}_{\tilde{E}}}Y = \lambda Y,
\end{equation}
therefore
$ \left(\operator{I}-\operator{I}_p\operator{Q}_{\tilde{E}}\right)Y = \lambda Y,$
and $\lambda$ is also an eigenvalue for the operator $\operator{M}=\operator{I}-\operator{I}_p\operator{Q}_{\tilde{E}}$.

To conclude, we will find a condition under which $M^k$ converges as  $k \rightarrow \infty$, which is a sufficient condition for 
 the eigenvalues of $\operator{M}$, other than one, to be strictly less than one in absolute value.  To do this, we define $\operator{\epsilon}=\operator{I}-\operator{I}_p$, and rewrite $\operator{M}$
as
\begin{align}
\operator{M}&=\operator{I}-\left(\operator{I}-\operator{\epsilon}\right)\qproje \text{,}
\end{align}
and 
\begin{align}
\operator{M}&=\proje+\operator{\epsilon}\qproje.
\end{align}
Since $\proje$ is a projection,
\begin{equation}
\qproje\proje=(\operator{I}-\proje)\proje=0.
\label{eq:qpzero}
\end{equation}
Squaring $\operator{M}$, we get
\begin{equation}
\operator{M}^2=\proje^2+\proje\operator{\epsilon}\qproje+\operator{\epsilon}\qproje\proje+\left(\operator{\epsilon}\qproje\right)^2,
\end{equation}
which, using the fact that $\proje^2=\proje$ and the orthogonality condition \eqref{eq:qpzero}, we simplify into
\begin{equation}
\operator{M}^2=\proje\left(\operator{I}+\operator{\epsilon}\qproje\right)+\left(\operator{\epsilon}\qproje\right)^2.
\end{equation}
By recurrence, we obtain that
\begin{equation}
\operator{M}^k=\proje\left(\sum_{i=0}^{k-1}\left(\operator{\epsilon}\qproje\right)^i\right)+\left(\operator{\epsilon}\qproje\right)^k,
\end{equation}
from which we conclude that the limit of $\operator{M}^k$  as $k\rightarrow\infty$ exists if and only if $\rho(\epsilon \qproje)$, the spectral radius of $\epsilon \qproje$, is less than 1, \ie
\begin{equation}
\rho\left(\epsilon \qproje \right)< 1.
\label{eq:stabilitycondition}
\end{equation}
\qed{}

Before we proceed to the numerical validations, we should discuss the implication of the above conditions. First, we should note that the above condition is satisfied if $\norm{\epsilon\qproje}<1$, and so for the condition \eqref{eq:stabilitycondition} to be satisfied, it is sufficient to have 
\begin{equation}
\norm{\epsilon}< \frac{1}{\norm{\qproje}}.
\label{eq:stab_cond2}
\end{equation}
Since the staggered operators $\proje$ and $\qproje$ on uniform grids are orthogonal with norm equal to one, we expect their norms on adaptive grids to be close to one (\ie $\norm{\qproje}\approx 1$), and so we expect the above stability condition to be met if $\norm{\epsilon}\lessapprox 1 $. 
In other words, we expect that if the norm of the interpolation error introduced by $\operator{I}_p$ is less than one, which is a realistic expectation as $\operator{I}_p$ is the product of consistent interpolations, and therefore consistent itself, then the projection is converging. 

Intuitively, one could think that for the iterated projection to converge, $\operator{I}_p$ must be stable and that, therefore, its norm must be less than one. Our proof suggests that this intuition is wrong and that, somewhat surprisingly, an unstable interpolation can lead to a converging algorithm as long as \eqref{eq:stab_cond2} is satisfied. This is essential as $\operator{I}_P$ involves interpolating from the nodes to the ghost nodes, which, as formula \eqref{eq:ghostnode_int} illustrates, is done using quadratic interpolating polynomials and can virtually introduce instabilities. The other two interpolations ($\operator{I}_{\tilde{N}}^{\tilde{E}}$ and $\operator{I}^{{N}}_{\tilde{E}}$) involved in $\operator{I}_P$ are constructed by taking weighted averages, and thus stable.

Our proof does not rely on a specific form for the ghost node interpolation; rather, it tells whether, for a given interpolation, the projection is stable. Unfortunately, our general ghost node construction is too complicated for us to prove whether condition \eqref{eq:stabilitycondition} is met or not. Even for simpler constructions, such as a linear one, where the ghost values are obtained as the weighted average of existing values, $\operator{I}_P$ is now the product of three averaging stochastic matrices, and thus itself a square stochastic matrix with positive coefficients, we can only conclude that its eigenvalues are in the disk $D(0,1)$, which is insufficient. 

Overall, the condition \eqref{eq:stabilitycondition} feels achievable and, in fact, less restrictive than one would have naively thought, and so we are reasonably hopeful that our method is stable. Yet, the constructions of our interpolation operators, even in the low-order scenario, are too complicated 
for us to prove it formally. In addition, we should point out that boundary effects are left out of our theoretical study. We, therefore, resort to an in-depth computational verification.

\subsection{Computational verification}
In this section, we verify the stability of our nodal projection operator by computing its convergence properties using a standard example with homogeneous Neumann boundary conditions and by verifying the stability of our nodal projection operator in the presence of different combinations of boundary and interface conditions. For the examples in this section and subsequent tests, wall and boundary conditions are treated using the second-order hybrid Finite-Volume/Finite Difference approach presented in \cite{THEILLARD2019108841}, itself based on previous studies \cite{Theillard2013AMM}. For interface conditions, the fluid quantities in $\Omega^-$ are systematically extended to $\Omega^+$ using the PDE-based approach proposed by Aslam in \cite{ASLAM2004349}. All other problem-specific treatments are described in the appropriate sections.

\subsubsection{Projection test - supra-convergence of the Hodge variable} 
We verify the convergence and stability of our projection operator by repeatedly applying the operator to an initially non-divergence-free velocity field and computing the error between the resulting velocity field and the exact divergence-free velocity field, as was done in previous studies \cite{min2006second,NG20098807,GUITTET2015215}. Doing so, we test the stability of our projection and the convergence of the Hodge variable as the grid resolution goes to zero. We consider the two-dimensional velocity field:
\begin{align}
\label{eq:non_div_free_vel}
    u^*(x,y)= &\sin(x)\cos(y)+x(\pi-x)y^2 \left (\frac{y}{3}-\frac{\pi}{2} \right ) \text{,} \\
    v^*(x,y)= & -\cos(x)\sin(y)+y(\pi-y)x^2 \left (\frac{x}{3}-\frac{\pi}{2} \right ) \text{,}
\end{align}
in the domain \(\Omega = [0,\pi]^2\) with Neumann boundary conditions. The divergence-free velocity field corresponding to \(u^*\) and \(v^*\) above is:
\begin{align}
\label{eq:div_free_vel}
    u(x,y)= & \sin(x)\cos(y) \text{,} \\
    v(x,y)= & -\cos(x)\sin(y).
\end{align}

\begin{figure}
    \centering
    \includegraphics[width=0.5\textwidth]{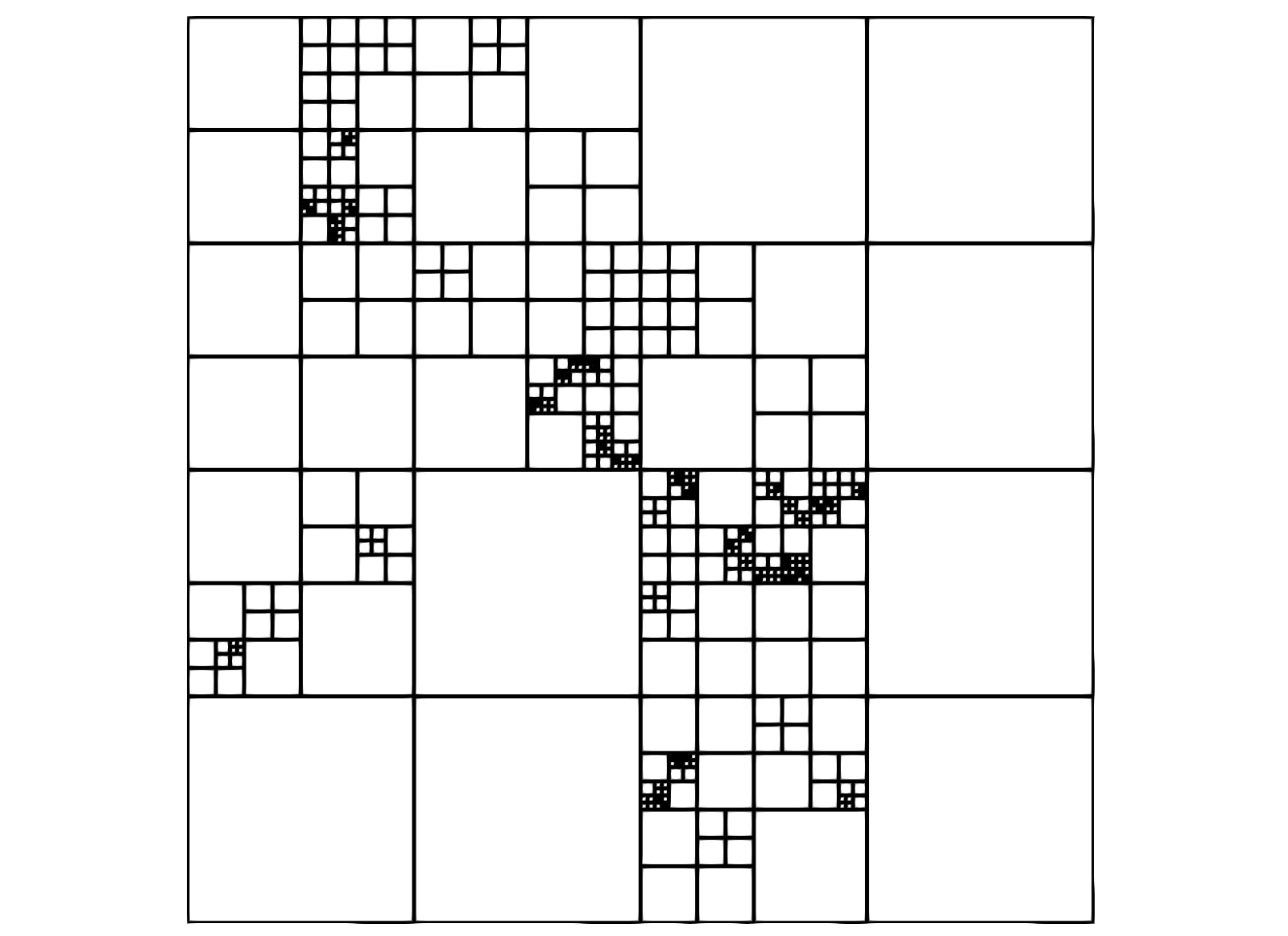}\\
    \captionof{figure}{Initial highly non-graded quadtree used for the Projection Test initialized with 240 random splits}
         \centering
         \includegraphics[width=0.45\textwidth]{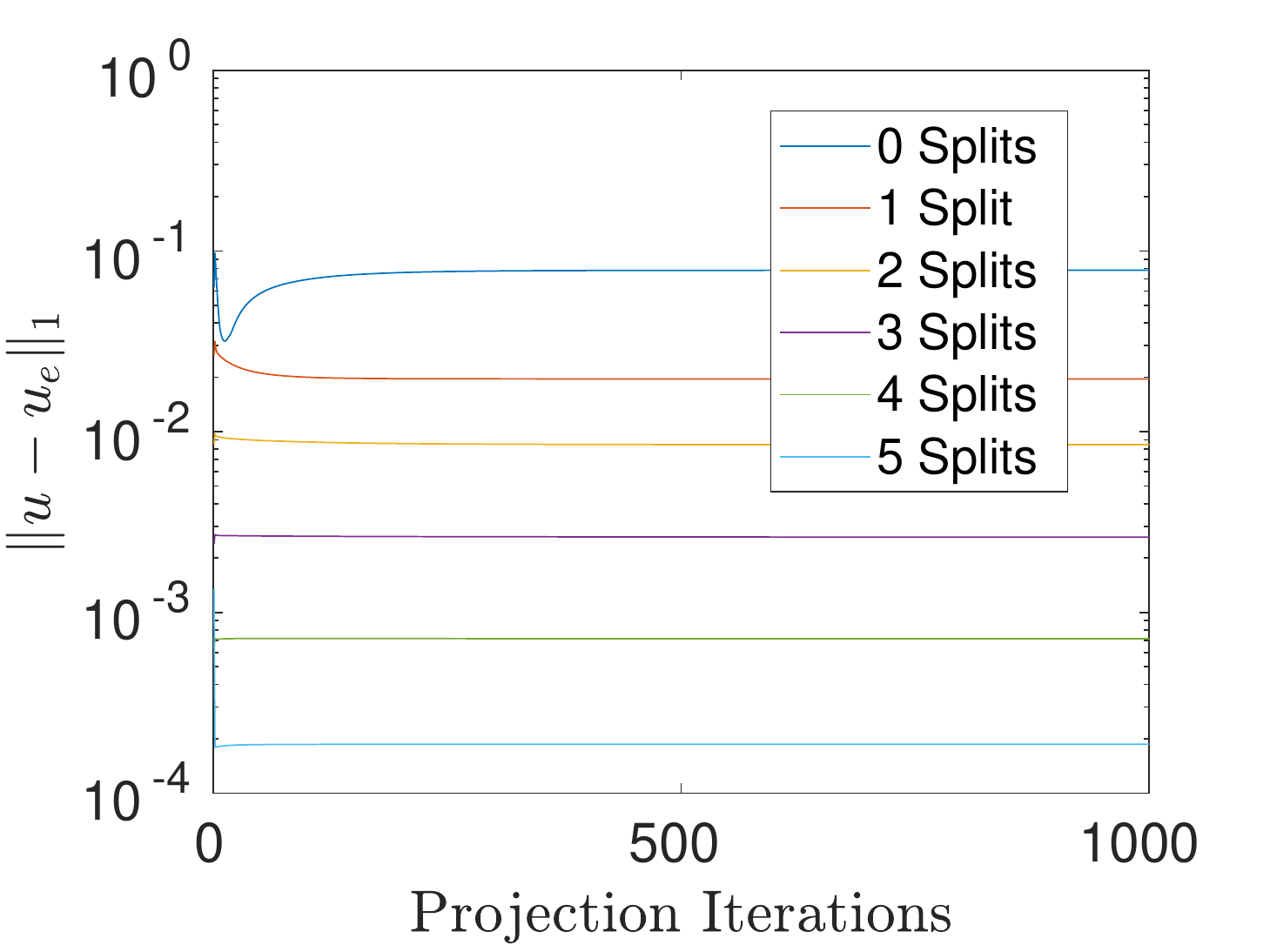}
         \centering
         \includegraphics[width=0.45\textwidth]{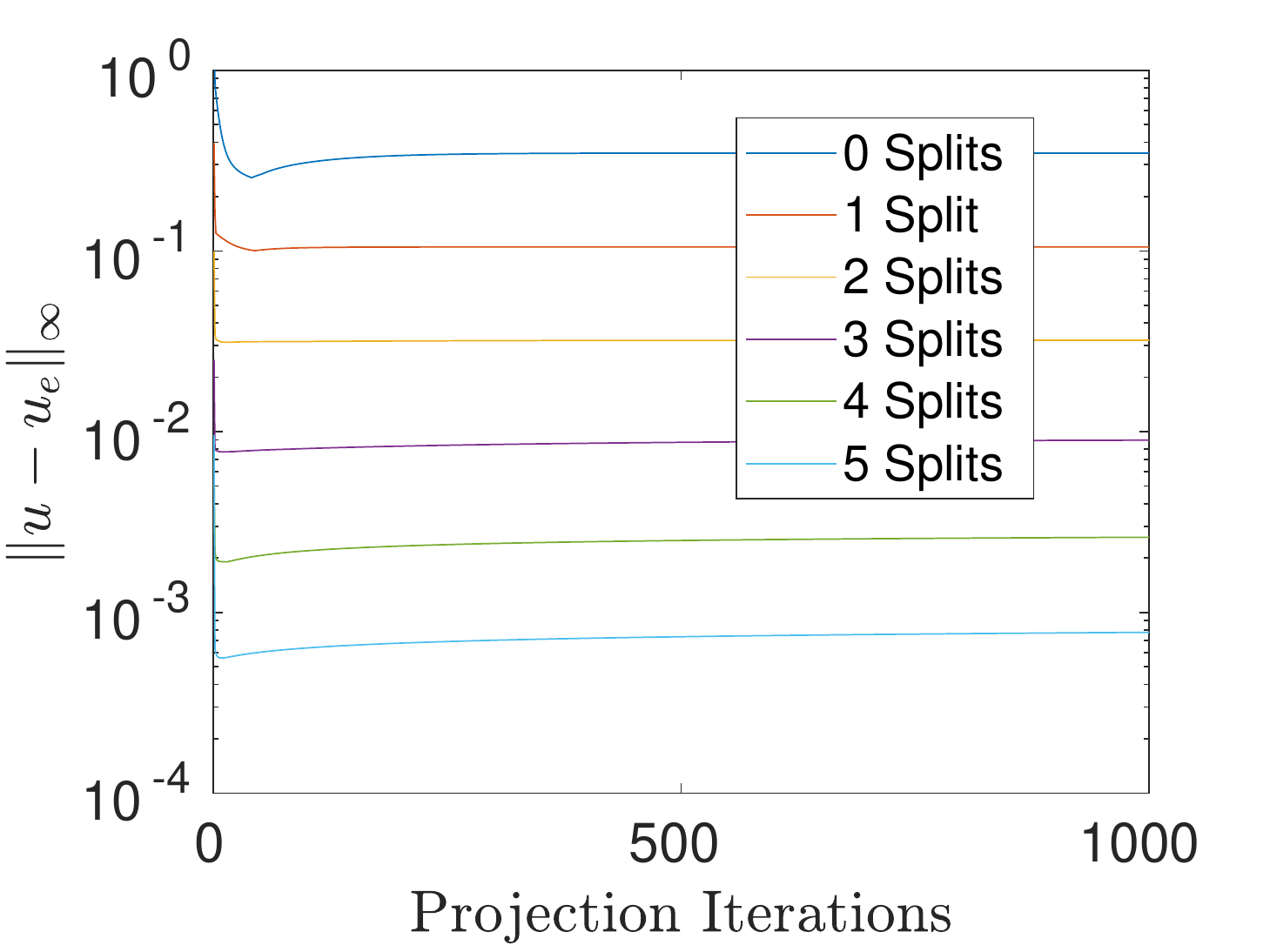}
         
    \caption{The $L^1$ (left) and $L^\infty$ (right) errors for the projection test on a randomly generated quadtree with increasing number of splits.} \label{fig:quadtree_splits}
    
    \captionof{table}{The final error and order of convergence for the projection test on a randomly generated quadtree with an increasing number of splits.}
        \begin{tabular}{|c|cccc|}
 
        \hline
        Splits      & $L^1$         & Order     & $L^\infty$     & Order    \\ \hline
        0           & 7.82e-02      & -         & 3.48e-01       & -        \\
        1           & 1.96e-02      & 2.00      & 1.05e-01       & 1.73     \\
        2           & 8.47e-03      & 1.21      & 3.20e-02       & 1.72     \\
        3           & 2.61e-03      & 1.70      & 8.99e-03       & 1.83     \\
        4           & 7.16e-04      & 1.87      & 2.61e-03       & 1.79     \\
        5           & 1.87e-04      & 1.94      & 7.76e-04       & 1.75     \\
        \hline
        \end{tabular}
        \label{fig:my_label}
    \label{fig:projection_test}
\end{figure}

We perform this test by starting with a randomized quadtree with 240 initial splits. Our goal here is to verify the order of accuracy for our projection operator in extreme scenarios where the grid is highly non-graded. We then successively refine the grid by further splitting each leaf, each time splitting the local spatial resolution in half, and monitor the error. The results of this test are shown in Figure \ref{fig:quadtree_splits}, along with the initial random quadtree. As expected, the projection is stable, and we see second-order convergence in the velocity field. 

\subsubsection{Projection test - stability} 
In this test, we specifically look at the stability of our projection operator in the presence of a variety of boundary and interface conditions. Using the same initial velocity field as in the previous example, \eqref{eq:non_div_free_vel}, we successively apply the projection operator to this field and monitor the norm of the variation between the velocity at the current projection iteration and the velocity after the final projection iteration. For a stable operator, this variation will tend to zero, up to the solver tolerance.

\begin{figure}[h!]
    \centering
     \includegraphics[width=0.62\textwidth]{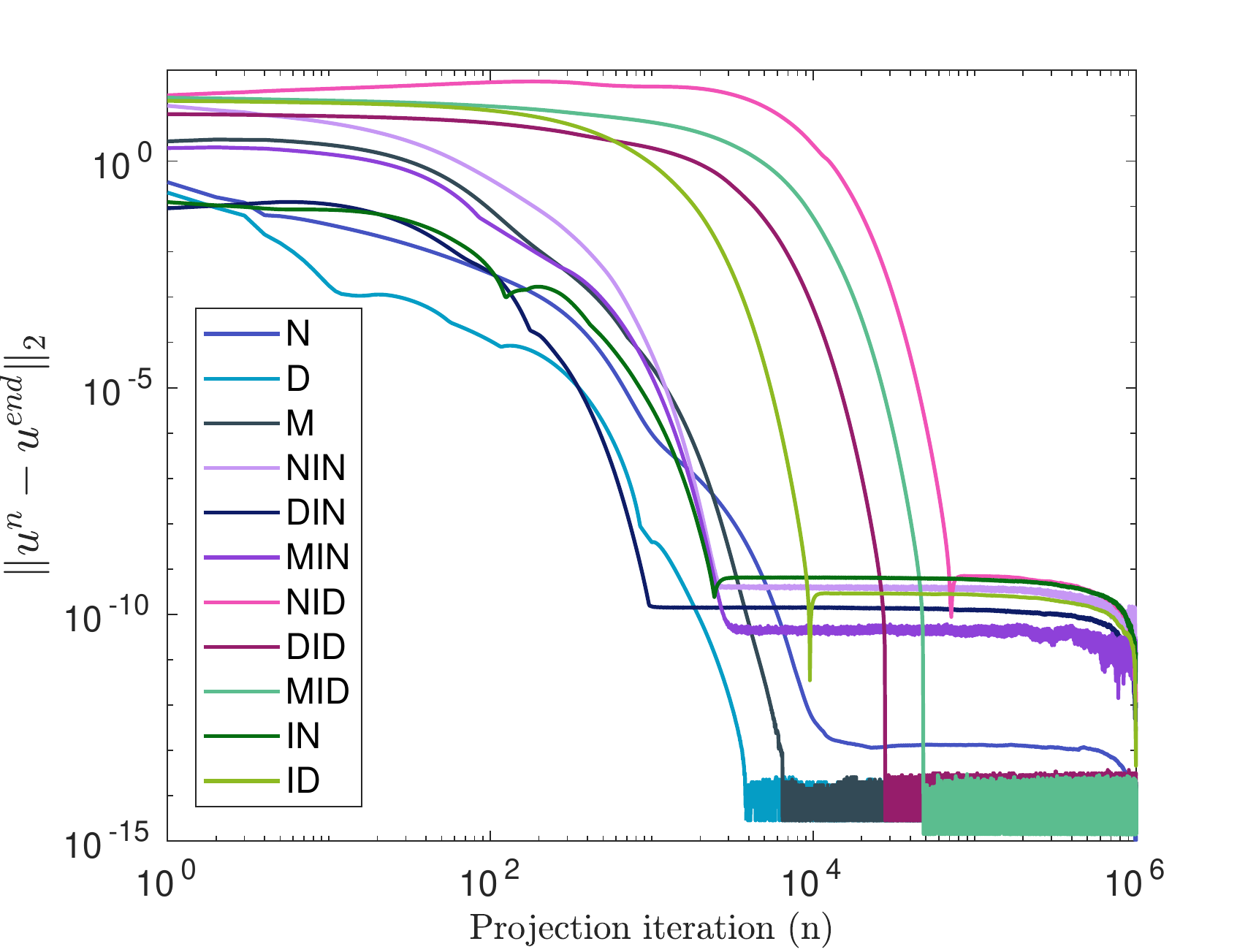}
     \includegraphics[width=0.92\textwidth]{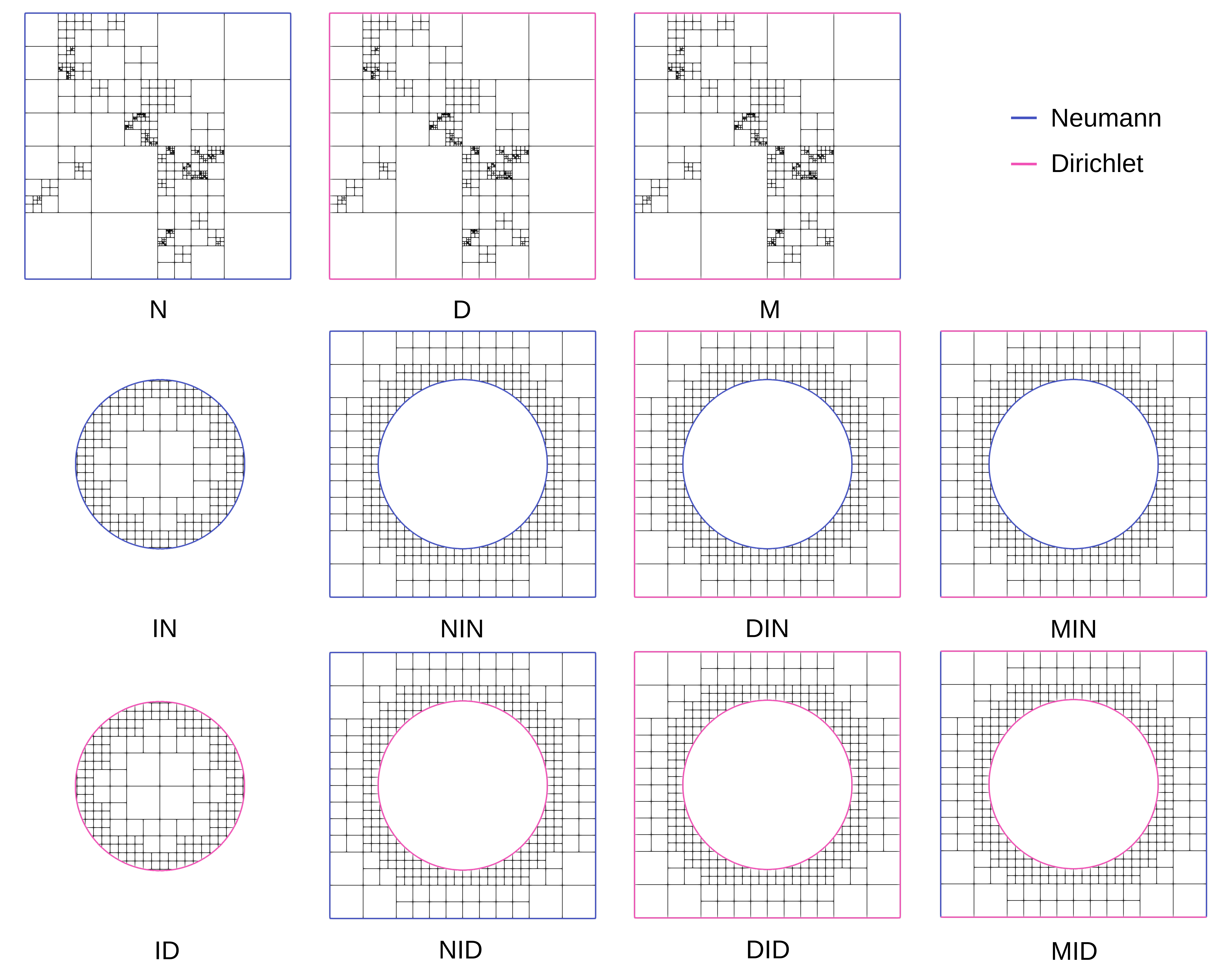}
    \caption{Variation in $L^2$ norm of the velocity after successive iterations of the projection operator for different sets of boundary conditions. Wall boundary conditions sets are denoted N for Neumann, D for Dirichlet, and M for mixed (a combination of Dirichlet and Neumann). When an interface is present, the wall boundary condition sets are augmented with IN for a Neumann interface and ID for a Dirichlet interface.}
    \label{fig:projection_stability}
\end{figure}

The results of this test are shown in Figure \ref{fig:projection_stability}. As expected, the variation in the norm of the velocity decreases as we successively apply our projection operator. For these examples, the tolerance of our linear solver was set to $10^{-12}$. As expected, we see that our collocated projection operator is numerically stable for all boundary and interface conditions tested. We also note that in practice, only a small number of iterations of the projection operator will be used and, as we later show, second-order accuracy in both the velocity and Hodge variable can be achieved with only a single projection applied.

\section{Nodal Navier-Stokes solver} \label{sec:NSsolver}
In this section, we integrate our nodal projection operator into a general framework for solving the incompressible Navier-Stokes equations. We begin by presenting an overview of our time-stepping algorithm and then we provide a brief introduction of the specific numerical techniques used to compute the intermediate velocity field. We conclude this section by solving an analytical solution to the Navier-Stokes equations with our nodal solver and comparing these results with the MAC-based solver presented in \cite{GUITTET2015215}. Furthermore, we use this example to discuss the practical number of projection iterations that need to be applied when using our nodal solver. 

\subsection{Overview of the time stepping algorithm}

We advance our velocity field from $\vectorr{u}^n$ to $\vectorr{u}\npo$ using a standard two-step projection method. In the Viscosity step, we compute the intermediate velocity field, $\mathbf{u^*}$, using a semi-Lagrangian Backward Difference Formula (SLBDF) scheme for the temporal integration of the momentum equation where the viscous terms are treated implicitly. In the Projection step, we repeatedly apply our nodal projection operator with the stopping criteria described in Section \ref{sec:corr_bc}. To account for the splitting error on the boundary (see Section \ref{sec:bcs}), we iterate the Viscosity and Projection steps following a boundary correction procedure detailed in Section \ref{sec:corr_bc}. A summary of our time-stepping algorithm can be seen in Figure \ref{fig:overall_algo} and details of the specific techniques used can be found in the subsequent sections.

\begin{figure}[!h]
\small
\begin{tabular}{c }
\hline \\ {\begin{minipage}[c]{0.9\textwidth}
\begin{description}
\item[1 -- Initialization] \hfill \\
Initialize the corrective velocity boundary condition $\vectorr{c}$ as the final one from the previous iteration.
\item[2 -- Repeat until $\norm{\vectorr{u}^{n+1}-\vectorr{u}|_\Gamma}<\epsilon_o$] \hfill 
\begin{description}
\item [2a -- Viscosity step] \hfill \\
Compute the intermediate velocity field $\mathbf{u}^*$ as the numerical solution of 
\begin{align}
\label{eq:visco_eq} \rho \frac{D {\vectorr{u}^*}}{Dt} &= \mu \triangle {\vectorr{u}^*}  +\vectorr{f} \qquad \forall \mathbf{\vectorr{x}} \in \Omega^-\slash \Gamma, \\
\label{eq::jump_u*}\vectorr{u}^* &= \vectorr{u}|_\Gamma + \vectorr{c} \qquad \forall \mathbf{\mathbf{x}} \in \Gamma,
\end{align}
Initialize $\vectorr{u}^{n+1}=\vectorr{u}^*$
\item [2b -- Projection step] \hfill
\begin{description}
    \item [ Repeat until $\norm{\vectorr{u}^{n+1} - \operator{P}_N \vectorr{u}^{n+1}}<\epsilon_i \norm{\vectorr{u}^{n+1}}$ ]\hfill \\
Project the velocity field
\begin{equation}
\vectorr{u}^{n+1} = \operator{P}_N \vectorr{u}^{n+1} \text{.}
\end{equation}
 \item [ Compute the new correction] \hfill \\ 
 \begin{equation}
 \label{eq:bc_corr}
     \vectorr{c}=\vectorr{c} - \omega \left(\vectorr{u}^{n+1} -\vectorr{u}|_{\Gamma}\right) \text{.}
 \end{equation}
\end{description}
\end{description}
\item [3 -- Update]\hfill \\
	Adapt the mesh to $\vectorr{u}^{n+1}$ and update all the variables accordingly.
\end{description}
\end{minipage}}
\\
\\
\\
\hline
\end{tabular}
	\caption{Outline of the algorithm for the construction of the solution $\mathbf{u}^{n+1}$ at time $t_{n+1}$ from the solution $\vectorr{u}^n$ at the previous time step $t_n$.}
	\label{fig:overall_algo}
	\end{figure}

\subsubsection{Viscosity step - SLBDF integrator }
\label{sec:sl_method}

To solve for the auxiliary velocity field $\vectorr{u}^*$, we start by integrating the momentum equation \eqref{eq:NS_mom} using a semi-Lagrangian backward difference (SLBDF) scheme \cite{long2013general,GUITTET2015215,THEILLARD201991} with an adaptive time-step $\Delta t_n = t_{n+1}-t_n$, and implicit treatment of the viscous term, namely
\begin{equation}
\rho \left( \alpha \frac{\vectorr{u}^*-\vectorr{u}^n_d}{\Delta t_n} +\beta \frac{\vectorr{u}^n_d - \vectorr{u}^{n-1}_d}{\Delta t_{n-1}}   \right) = \mu \triangle \vectorr{u}^* +\vectorr{f},
\end{equation}
where 
\begin{equation}
    \alpha = \frac{2\Delta t_n + \Delta t_{n-1}}{\Delta t_n +\Delta t_{n-1}} \text{,} \ \ \beta = -\frac{\Delta t_n}{\Delta t_n + \Delta t_{n-1}}.
\end{equation}
The departing velocities $\vectorr{u}^n_d$ and $\vectorr{u}^{n-1}_d$ are obtained by interpolating the velocity field at the departure points $\vectorr{x}_{d}^{n}$ and $\vectorr{x}_{d}^{n-1}$ and corresponding time steps. The adaptive time-step, $\Delta t_n$ is determined by pre-setting the CFL number, $\textrm{CFL} = \max\norm{\vectorr{u}}\Delta t/\Delta x$. To find the departure points, we follow the characteristic curve backward using a second-order Runge-Kutta method (RK2, midpoint): 
\begin{align}
    \hat{\vectorr{x}} &= \vectorr{x}^{n+1} - \frac{\Delta t_n}{2}  \cdot \vectorr{u}^{n+1}(\vectorr{x}^{n+1}), \label{eq:xnd_intermediate} \\
    \vectorr{x}^n_d &= \vectorr{x}^{n+1} - \Delta t_n  \cdot\vectorr{u}^{n+\frac{1}{2}}(\hat{\vectorr{x}}),\label{eq:xnd_departure}
\end{align}
and
\begin{align}
    \Bar{\vectorr{x}} &= \vectorr{x}^{n+1} - \frac{(\Delta t_n + \Delta t_{n+1})}{2} \cdot \vectorr{u}^{n+1}(\vectorr{x}^{n+1}), \label{eq:xnm1d_intermediate} \\
    \vectorr{x}^{n-1}_d &= \vectorr{x}^{n+1} - (\Delta t_n + \Delta t_{n+1}) \cdot \vectorr{u}^{\textrm{mid}}(\Bar{\vectorr{x}}).\label{eq:xnm1d_departure}
\end{align}
 We interpolate the intermediate velocities, $\vectorr{u}^{n+\frac{1}{2}}(\hat{\vectorr{x}})$ and $\vectorr{u}^{\textrm{mid}}(\Bar{\vectorr{x}})$, from the velocity fields at $t_{n-1}$ and $t_n$ as
\begin{align}
    \vectorr{u}^{n+\frac{1}{2}}(\hat{\vectorr{x}}) &= \frac{2 \Delta t_{n-1} + \Delta t_{n}}{2 \Delta t_{n-1}} \vectorr{u}^{n}(\hat{\vectorr{x}}) - \frac{\Delta t_n}{2 \Delta t_{n-1}} \vectorr{u}^{n-1}(\hat{\vectorr{x}}) \text{,}\\
    \vectorr{u}^{\text{mid}}(\bar{\vectorr{x}}) &= \frac{\Delta t_n + \Delta t_{n-1}}{2\Delta t_{n-1}}\vectorr{u}^n(\bar{\vectorr{x}}) + \frac{\Delta t_{n-1}-\Delta t_n}{2\Delta t_{n-1}}\vectorr{u}^{n-1}(\bar{\vectorr{x}}) \text{.}
\end{align}
\subsubsection{Viscosity step - improved trajectory reconstruction} 
\label{sec:improved_sl}

 As the semi-Lagrangian method shown in Section \ref{sec:sl_method} relies on knowing what $\vectorr{u}^{n+1}$ is, we must look to other explicit methods to find the departure point. Here, we discuss two departure point methods that utilize an RK2 scheme and approximate $\vectorr{u}^{n+1}.$
 
 The first and most obvious choice for utilizing an RK2 method with known information is to replace $\vectorr{u}^{n+1}$ in \eqref{eq:xnd_intermediate} and \eqref{eq:xnm1d_intermediate} with $\vectorr{u}^n.$ This midpoint rule formulation was popularized by Xiu and Karniadakis \cite{XIU2001658} and has been a popular integration choice for projection methods that utilize a semi-Lagrangian discretization of the advection term (see, \eg \cite{GUITTET2015215,min2006second}). Using $\vectorr{u}^n$ is akin to a constant Taylor approximation to $\vectorr{u}^{n+1}.$ 
 
 This constant approximation is often a sufficient approximation for lower CFL numbers, but for larger CFL numbers, it can cause the computed departure point to lie too far from the characteristic prior to interpolation and can lead to less accurate velocity terms.  By doing a first-order Taylor expansion of $\vectorr{u}^{n+1}(\vectorr{x}^{n+1})$ around $t_n$ and using an Euler approximation for any derivative terms, we get
\begin{equation}
    \vectorr{u}^{n+1}(\vectorr{x}^{n+1})=\vectorr{u}^n(\vectorr{x}^{n+1}) + \frac{\Delta t_n}{\Delta t_{n-1}} \left(\vectorr{u}^{n}(\vectorr{x}^{n+1})-\vectorr{u}^{n-1}(\vectorr{x}^{n+1})\right) + \mathcal{O}(\Delta t^2),
\end{equation}
which we substitute into \eqref{eq:xnd_intermediate} and \eqref{eq:xnm1d_intermediate}. By having a higher-order approximation for $\vectorr{u}^{n+1}(\vectorr{x}^{n+1}),$ we expect our numerical simulations to have more accurate results, especially if a large time-step is used. 

 While we can mitigate this issue with a higher-order integration scheme, such as a fourth-order Runge-Kutta method, the order of accuracy of the semi-Lagrangian method is determined by both the integration scheme used to find the departure point and the interpolation scheme used to evaluate $\vectorr{u}$ at the departure point. Since we are using a quadratic interpolation scheme, using an integration scheme that is higher than second-order accurate is superfluous.

We use both integrating schemes for our numerical validation in Section \ref{sec:examples}. The classical scheme is used to directly compare the nodal projection method with previous numerical studies and experiments, and the new improved integrating scheme is used to compare the accuracy of the solver when implementing each integrating scheme.

\subsubsection{Stopping criteria for iterative procedures} \label{sec:corr_bc}
As our projection operator is not an exact orthogonal projection, we use an iterative procedure to create an approximately divergence-free velocity field. We perform successive projections until 
\begin{equation}
\norm{\vectorr{u}^{n+1} - \operator{P}_N \vectorr{u}^{n+1}}<\epsilon_i \norm{\vectorr{u}^{n+1}} \text{,}
\end{equation}
or a predefined maximum number of iterations, $K_{max}$, has been reached. Typically, we choose $\epsilon_i=10^{-3}$, set $K_{max} = 5$, and only a few ($1-3$) iterations are required to reach convergence. 

The boundary correction $\vectorr{c}$ is designed to compensate for the splitting errors (see \ref{sec:bcs}) in a similar fashion to what was proposed in \cite{THEILLARD201991}. It is designed so that when the correction reaches convergence (see Eq. (\ref{eq:bc_corr})), the boundary condition on the solid object is satisfied (\ie $\vectorr{u}^{n+1} = \vectorr{u}$). The parameter $\omega$ controls the convergence rate and must be chosen in the range $0<\omega<1$. Similar corrections are performed on the wall of the computational domain $\partial \Omega$ but are not explicited here for the sake of concision. Typically, we terminate the outer iterations when the error in the interface's velocity ($\norm{\vectorr{u}^{n+1}-\vectorr{u}|_\Gamma}$) is less than $\epsilon_o=10^{-3}$. 

\subsubsection{Refinement criteria}
\label{sec:refinement}
As it was done in our previous studies \cite{GUITTET2015215,THEILLARD201991,THEILLARD2019108841,CLERETDELANGAVANT2017271,THEILLARD2021110478}, the mesh is dynamically refined near the solid interface and where high gradients of vorticity occur. At each iteration, we recursively apply the following splitting criterion at each cell. 

We split each cell $\mathcal{C}$ if
\begin{eqnarray}
    \min_{n \in \text{nodes} (\mathcal{C})} \left | \phi(n) \right | \leq B\cdot\text{Lip}(\phi) \cdot \text{diag}(\mathcal{C}) 
    \quad \text{and} \quad 
    \text{level}(\mathcal{C}) \leq \textrm{max}_{\textrm{level}} \text{,}
\end{eqnarray}
or 
\begin{eqnarray}
    \min_{n \in \text{nodes} (\mathcal{C})} \text{diag}(\mathcal{C})\cdot\frac{\norm{\nabla\vectorr{u}(n)}}{\norm{\vectorr{u}}_\infty} \geq T_V 
    \quad \text{and} \quad 
    \text{level}(\mathcal{C}) \leq \textrm{max}_{V} \text{,}     
\end{eqnarray}
or
\begin{eqnarray}
    \text{level}(\mathcal{C}) \leq \textrm{min}_{\textrm{level}}.\label{eq:refcritminlevel}
\end{eqnarray}
If none of these criteria are met, we merge $\mathcal{C}$ by removing all its descendants. 

Here $\text{Lip}(\phi)$ is an upper estimate of the minimal Lipschitz constant of the level set function $\phi$. Since the level set used to create the mesh will be reinitialized (\ie $|\nabla \phi|=1$), we use  $\text{Lip}(\phi)=1.2$.  $\text{diag}(\mathcal{C})$ is the length of the diagonal of cell $\mathcal{C}$. $B$ is the user-specified width of the uniform band around the interface.  $T_V$ is the vorticity threshold and $\max_{\textrm{level}}$ is the maximum grid level allowed for the vorticity-based refinement. Condition (\ref{eq:refcritminlevel}) ensures that a minimum resolution of $\min_{\textrm{level}}$ is maintained.

\subsection{Verification - analytic vortex} \label{sec:avortex}
\subsubsection{Convergence study} 
\begin{figure*}[ht]
   
\centering

\begin{scriptsize}
\captionof{table}{Convergence of x-component velocity and Hodge variable $\Phi$ for the Nodes and MAC \cite{GUITTET2015215} implementations.}
\label{table:avortex}
\begin{tabular}{|c|cccc|cccc|}
\multicolumn{9}{c}{} \\
\multicolumn{9}{c}{\bf Velocity} \\ \hline
                            & \multicolumn{4}{c|}{Nodes}                            & \multicolumn{4}{c|}{MAC} \\
Level (max:min)    & $L^1$         & Order     & $L^\infty$    & Order     & $L^1$         & Order     & $L^\infty$    & Order     \\ \hline
7:3                         & 6.84e-04      & -         & 2.34e-02      & -         & 5.48e-03      & -         & 2.18e-02      & -         \\
8:4                         & 1.25e-04      & 2.46      & 2.72e-03      & 3.11      & 1.76e-03      & 1.64      & 8.65e-03      & 1.33      \\
9:5                         & 2.33e-05      & 2.42      & 9.56e-04      & 1.51      & 5.56e-04      & 1.66      & 3.28e-03      & 1.40      \\
10:6                        & 5.61e-06      & 2.05      & 2.19e-04      & 2.12      & 1.61e-04      & 1.79      & 1.64e-03      & 1.00      \\
11:7                        & 1.36e-06      & 2.05      & 4.44e-05      & 2.30      & 3.62e-05      & 2.15      & 3.84e-04      & 2.09      \\
12:8                        & 4.15e-07      & 1.71      & 9.55e-06      & 2.22      & 8.59e-06      & 2.08      & 2.45e-04      & 0.65      \\
\hline
\multicolumn{9}{c}{} \\

\multicolumn{9}{c}{\bf Hodge} \\ \hline
                            & \multicolumn{4}{c|}{Nodes}                         & \multicolumn{4}{c|}{MAC} \\
Level (max:min)    & $L^1$         & Order     & $L^\infty$    & Order     & $L^1$     & Order         & $L^\infty$    & Order     \\ \hline
7:3                         & 4.43e-05      & -         & 4.16e-04      & -         & 3.17e-04      & -         & 1.62e-03      & -         \\
8:4                         & 1.67e-05      & 1.41      & 9.48e-05      & 2.13      & 3.47e-05      & 3.19      & 2.47e-04      & 2.71      \\
9:5                         & 2.55e-06      & 2.71      & 1.44e-05      & 2.72      & 4.48e-06      & 2.95      & 5.85e-05      & 2.08      \\
10:6                        & 3.07e-07      & 3.05      & 2.45e-06      & 2.55      & 7.20e-07      & 2.64      & 1.44e-05      & 2.02      \\
11:7                        & 3.12e-08      & 3.30      & 5.49e-07      & 2.16      & 7.41e-08      & 3.28      & 2.81e-06      & 2.36      \\
12:8                        & 3.67e-09      & 3.09      & 1.05e-07      & 2.39      & 1.14e-08      & 2.70      & 6.72e-07      & 2.06      \\\hline
\end{tabular}

\end{scriptsize}
 \centering
         \centering
         \includegraphics[width=0.47\textwidth]{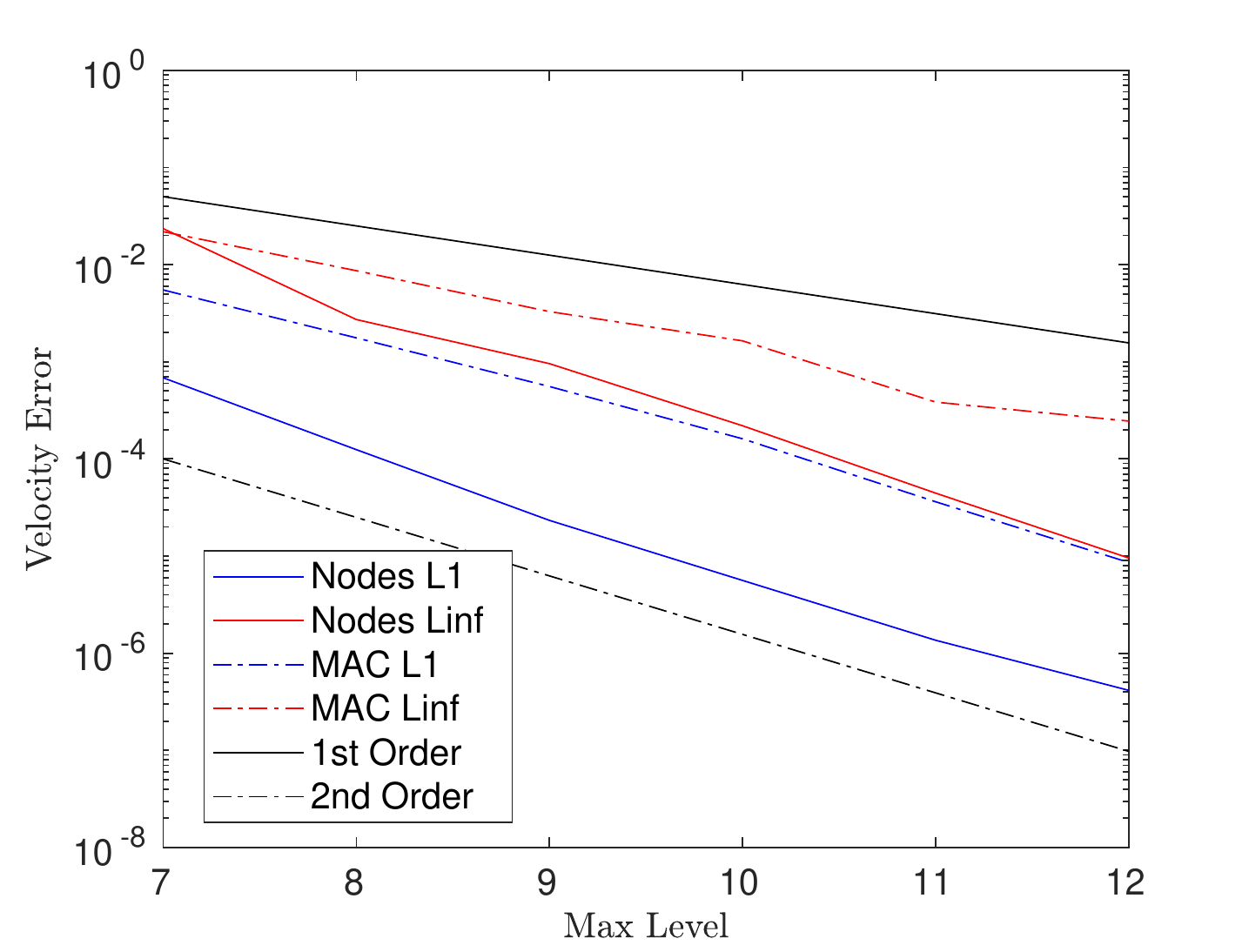}
     \hfill
         \centering
         \includegraphics[width=0.47\textwidth]{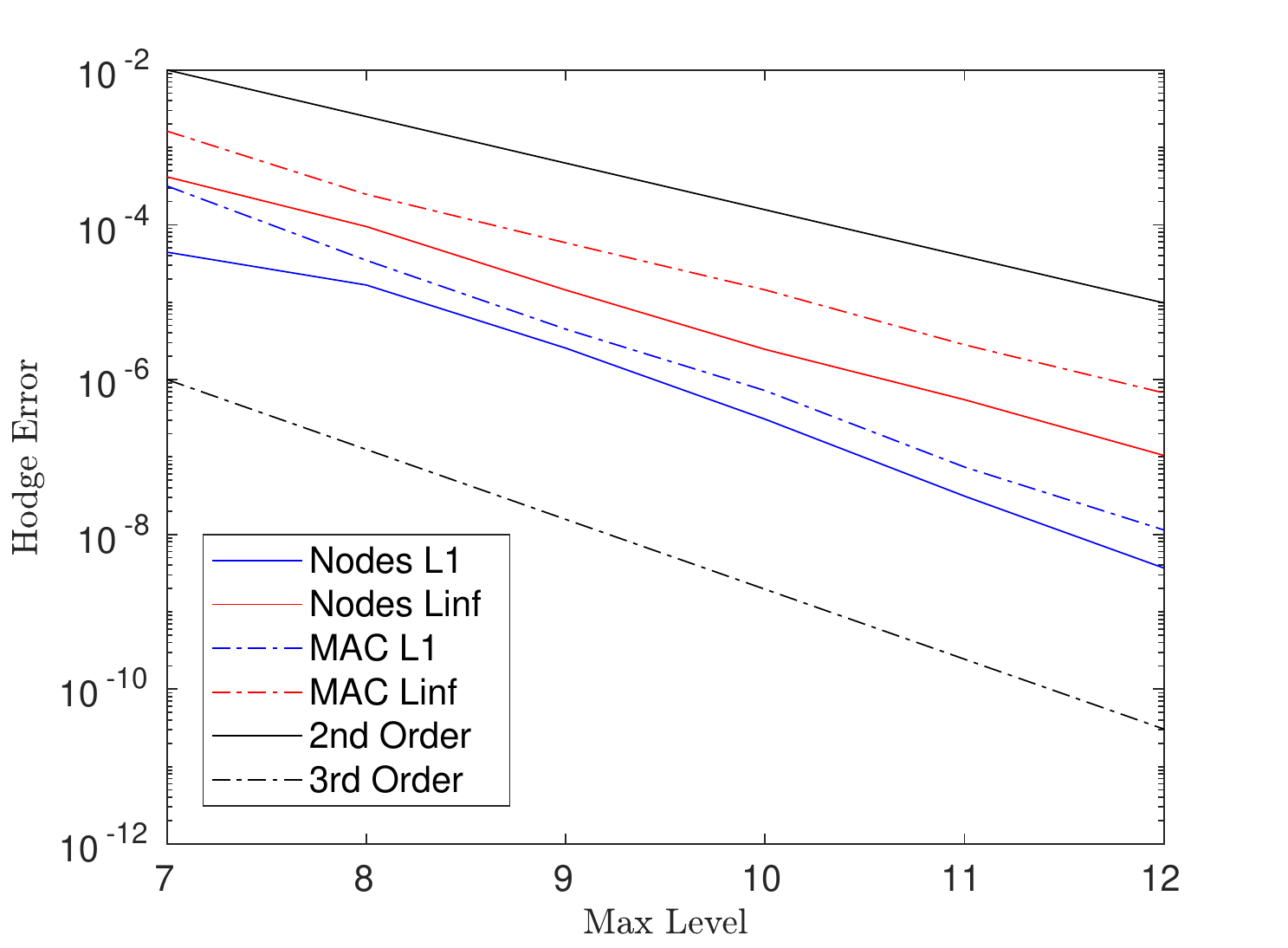}
    \caption{Visualization of the error convergence for the velocity and Hodge variable.}
    \label{fig:analytic_vortex}
\end{figure*}

We consider the following analytical solution to the Navier-Stokes Equations,
\begin{align}
    u(x,y) &=  \sin(x) \cos(y) \cos(t) \text{,} \\ 
    v(x,y) &=  -\cos(x) \sin(y) \cos(t) \text{,} \\
    p(x,y) &=  0 \text{.}
\end{align}
We set $\mu=1$, $\rho=1$, and define the forcing terms as
\begin{align}
    f_x &= \sin(x) \cos(y) \left (2 \mu \cos(t) - \rho \sin(t) \right ) + \rho \cos^2(t) \sin(x) \cos(x) \text{,} \\ 
    f_y &= \cos(x) \sin(y) \left (\rho \sin(t) - 2 \mu \cos(t) \right ) + \rho \cos^2(t) \sin(y) \cos(y).
\end{align}
We monitor the convergence of our solver as our mesh is refined and compare the results against the MAC-based solver of \cite{GUITTET2015215}. Each of the simulations is run until a final time of $t_f = \pi / 3$ is reached with an adaptive time step defined by $\Delta t_n = \Delta x / \max \| u_n \|$. We set the maximum number of projections, $K_{max}$, to 5, with an error tolerance of $10^{-3}$ (see Section \ref{sec:corr_bc}), and use the boundary correction procedure explained in Section \ref{sec:corr_bc}. The results of this example are shown in Table \ref{table:avortex} and in Figure \ref{fig:analytic_vortex}.

We see that our nodal projection method achieves second-order convergence in both the $L^1$ and $L^\infty$ norm for velocity as expected based on our discretization. For the Hodge variable, we also see second-order convergence in both $L^1$ and $L^\infty$. If we compare the results of our nodal solver with the MAC-based solver of \cite{GUITTET2015215}, we notice a distinct improvement in the $L^\infty$ norm for velocity. 

\subsubsection{Impact of successive projections - determining $K_{max}$}
\label{sec:det_kmax}
In order to develop a practical guideline for the number of projection iterations that should be used in our nodal solver, we again consider the analytic vortex example from the previous section. We continue using the same adaptive time step, grid refinement criteria, and boundary correction procedures previously described. However, instead of using the error threshold described in \ref{sec:corr_bc}, we explicitly define the number of projections to be performed at each time step. The results of this test are shown in Table \ref{tab:det_kmin} for the explicit number of projections set to 1, 5, and 10. 

We emphasize that for this example, only a single projection is required to achieve second-order accuracy in the $L^\infty$ norm for both velocity and the Hodge variable. Increasing the number of projections does not significantly change the solution or the order of accuracy. For this reason, we do not set a minimum number of projections for our nodal solver. Instead, we specify a maximum number of projections, $K_{max}$, and typically use a default value of 5 projections. In practice and using an error threshold of $10^{-3}$ (see Section \ref{sec:corr_bc}), we observe that a single projection is often sufficient.

\begin{figure*}[ht]
\centering
\begin{scriptsize}
\captionof{table}{Impact of the number of projections on the $L^\infty$ error for the velocity and Hodge variable $\Phi$ for the analytic vortex example.}
\begin{tabular}{|c|cccc|cccc|cccc|} \hline
         & \multicolumn{4}{c|}{Projections = 1}                             &  \multicolumn{4}{c|}{Projections = 5}                             & \multicolumn{4}{c|}{Projections = 10} \\
         & \multicolumn{2}{c}{Velocity}     & \multicolumn{2}{c|}{Hodge}    & \multicolumn{2}{c}{Velocity}      & \multicolumn{2}{c|}{Hodge}    & \multicolumn{2}{c}{Velocity}     & \multicolumn{2}{c|}{Hodge}     \\
Level    & $L^\infty$   & Order             & $L^\infty$    & Order         & $L^\infty$   & Order              & $L^\infty$    & Order         & $L^\infty$   & Order              & $L^\infty$    & Order         \\ \hline
7:3      & 2.36e-02     & -                 & 4.16e-04      & -             & 4.19e-02     & -                  & 9.65e-05      & -             & 5.44e-02     & -                  & 7.01e-05      & -             \\
8:4      & 2.79e-03     & 3.08              & 9.49e-05      & 2.13          & 6.86e-03     & 2.61               & 3.24e-05      & 1.58          & 8.74e-03     & 2.64               & 1.44e-05      & 2.28          \\
9:5      & 9.74e-04     & 1.52              & 1.44e-05      & 2.72          & 2.32e-03     & 1.56               & 5.81e-06      & 2.48          & 3.03e-03     & 1.53               & 2.89e-06      & 2.32          \\
10:6     & 2.25e-04     & 2.12              & 2.45e-06      & 2.56          & 4.92e-04     & 2.24               & 9.52e-07      & 2.61          & 6.38e-04     & 2.25               & 6.31e-07      & 2.20          \\
11:7     & 4.58e-05     & 2.29              & 5.49e-07      & 2.16          & 9.92e-05     & 2.31               & 2.11e-07      & 2.17          & 4.58e-05     & 3.80               & 5.48e-07      & 0.20          \\
12:8     & 9.96e-06     & 2.20              & 1.05e-07      & 2.38          & 1.82e-05     & 2.45               & 4.24e-08      & 2.32          & 2.45e-05     & 0.90               & 3.47e-08      & 3.98          \\
\hline
\end{tabular}
\label{tab:det_kmin}
\end{scriptsize}
\end{figure*}

\subsubsection{Verification of the improved departure point reconstruction scheme}

Finally, we use the analytic vortex example to verify that our improved departure point reconstruction scheme remains second-order accurate. Using an identical setup as in Section \ref{sec:avortex}, we monitor the convergence of our solver using the original SLBDF scheme presented in \cite{GUITTET2015215} and compare the results against the improved scheme presented herein. Again, the maximum number of projections, $K_{max}$, is set to 5, with an error tolerance for the projection set to $10^{-3}$. The $L^1$ and $L^\infty$ errors for the x-component velocity and the Hodge variable are shown in Table \ref{table:avortex_SLBDF}. As expected, we that the improved reconstruction scheme maintains second-order accuracy. 

\begin{figure*}[ ht]
\centering
\begin{scriptsize}
\captionof{table}{Convergence of the x-component velocity and Hodge variable $\Phi$ using the original and improved departure point reconstruction schemes.}
\label{table:avortex_SLBDF}
\begin{tabular}{|c|cccc|cccc|}
\multicolumn{9}{c}{} \\
\multicolumn{9}{c}{\bf Velocity} \\ \hline
                            & \multicolumn{4}{c|}{Original Reconstruction}                            & \multicolumn{4}{c|}{Improved Reconstruction} \\
Level (max:min)    & $L^1$         & Order     & $L^\infty$    & Order     & $L^1$         & Order     & $L^\infty$    & Order     \\ \hline
7:3                         & 6.84e-04      & -         & 2.34e-02      & -         & 9.08e-04      & -         & 2.90e-02      & -         \\
8:4                         & 1.25e-04      & 2.46      & 2.72e-03      & 3.11      & 1.13e-04      & 3.0.      & 2.73e-03      & 3.41      \\
9:5                         & 2.33e-05      & 2.42      & 9.56e-04      & 1.51      & 2.38e-05      & 2.24      & 9.58e-04      & 1.51      \\
10:6                        & 5.61e-06      & 2.05      & 2.19e-04      & 2.12      & 5.58e-06      & 2.09      & 2.20e-04      & 2.12      \\
11:7                        & 1.36e-06      & 2.05      & 4.44e-05      & 2.30      & 1.30e-06      & 2.10      & 4.47e-05      & 2.30      \\
12:8                        & 4.15e-07      & 1.71      & 9.55e-06      & 2.22      & 3.95e-07      & 1.72      & 9.55e-06      & 2.23      \\
\hline
\multicolumn{9}{c}{} \\

\multicolumn{9}{c}{\bf Hodge} \\ \hline
                            & \multicolumn{4}{c|}{Original Reconstruction}                         & \multicolumn{4}{c|}{Improved Reconstruction} \\
Level (max:min)    & $L^1$         & Order     & $L^\infty$    & Order     & $L^1$     & Order         & $L^\infty$    & Order     \\ \hline
7:3                         & 4.43e-05      & -         & 4.16e-04      & -         & 3.59e-05      & -         & 3.18e-04      & -         \\
8:4                         & 1.67e-05      & 1.41      & 9.48e-05      & 2.13      & 1.69e-05      & 1.09      & 9.72e-05      & 1.71      \\
9:5                         & 2.55e-06      & 2.71      & 1.44e-05      & 2.72      & 2.57e-06      & 2.72      & 1.47e-05      & 2.72      \\
10:6                        & 3.07e-07      & 3.05      & 2.45e-06      & 2.55      & 3.10e-07      & 3.05      & 2.44e-06      & 2.59      \\
11:7                        & 3.12e-08      & 3.30      & 5.49e-07      & 2.16      & 3.20e-08      & 3.28      & 5.59e-07      & 2.13      \\
12:8                        & 3.67e-09      & 3.09      & 1.05e-07      & 2.39      & 3.83e-09      & 3.06      & 1.08e-07      & 2.37      \\\hline
\end{tabular}
\end{scriptsize}
\end{figure*}

\section{Validation examples} \label{sec:examples}
In this section, we validate our node-based Navier-Stokes solver using several canonical problems in two and three spatial dimensions. We compare the results of our solver with previous numerical studies \cite{ghia1982high, erturk2005numerical, GUITTET2015215, NG20098807, braza_chassaing_minh_1986, CALHOUN2002231, ENGELMAN1990, MITTAL20084825, MARELLA20051, johnson_patel_1999, EGAN2021110084} and experimental results \cite{dutsch1998low} demonstrating the robustness of our solver. Furthermore, we use the Kármán vortex example to highlight the advantages of our improved departure point reconstruction scheme and the efficiency of our solver. Finally, we conclude this section by simulating a high Reynolds number flow past the UC Merced Beginnings sculpture to demonstrate that our solver can be used to simulate flows with non-trivial geometries. 

\subsection{Driven cavity} \label{sec:driven_cavity}
The lid-driven cavity problem is a canonical example for validating an incompressible flow solver. Following the numerical experiments of Ghia \etal \cite{ghia1982high} and Erturk \etal \cite{erturk2005numerical}, we consider the domain $\Omega = [0,1]^2$ with the top wall moving at a unit velocity and no-slip boundary conditions on the velocity for all four walls. We set the density $\rho = 1$ and the Reynolds number to $Re=1~000$.

\begin{figure}
         \includegraphics[width=0.47\textwidth]{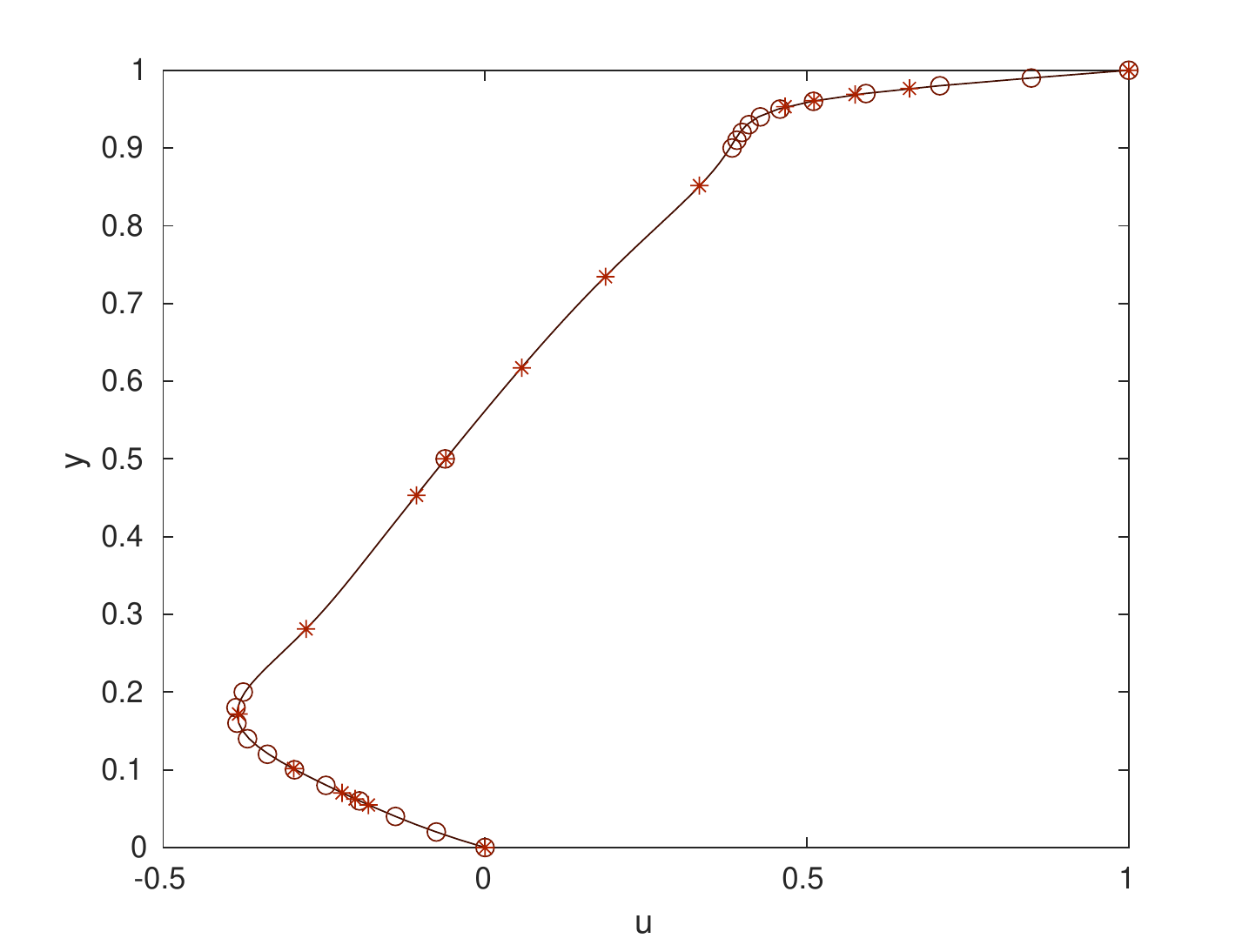}
         \includegraphics[width=0.47\textwidth]{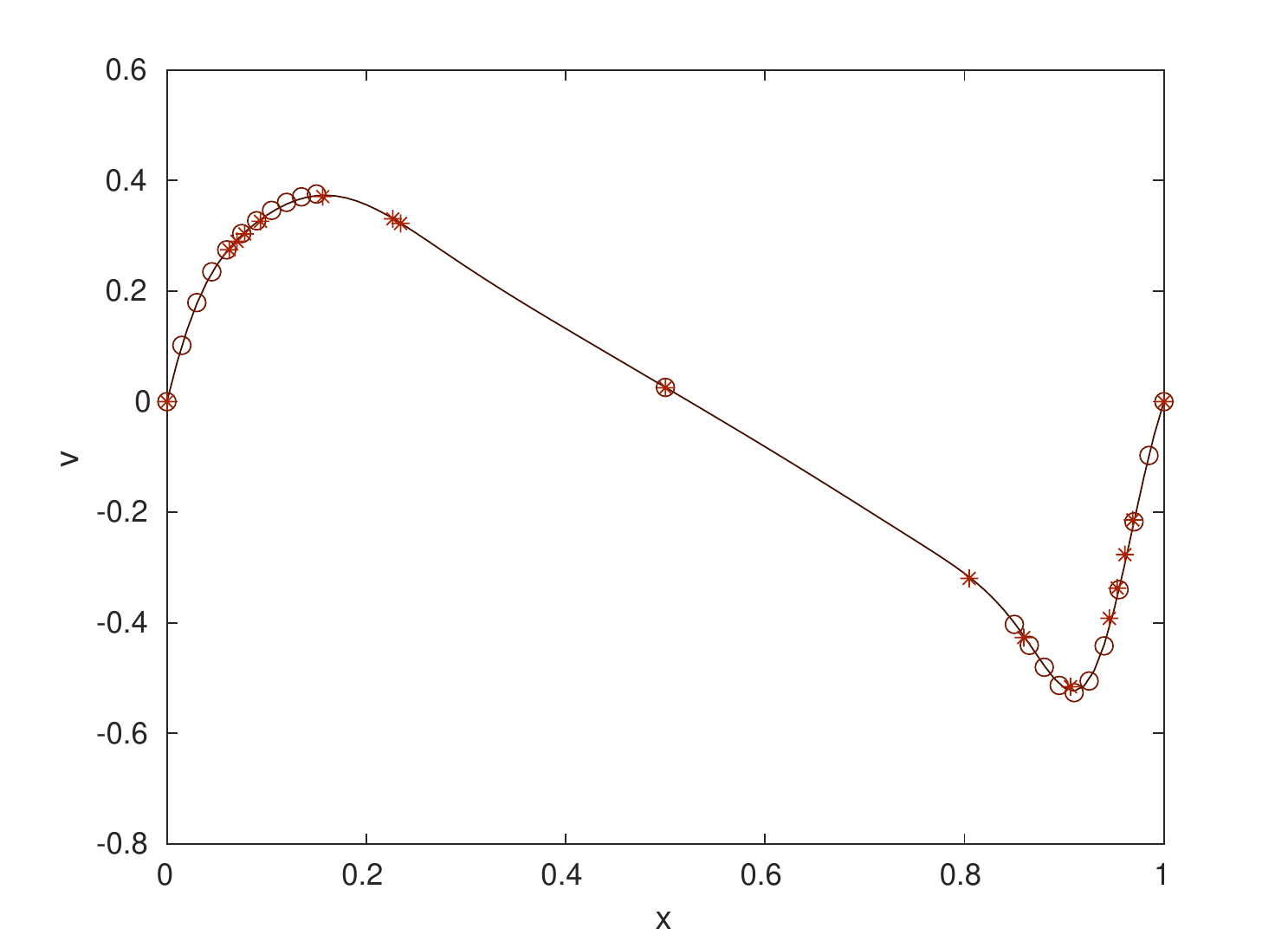}
    \caption{x- and y-components of the velocity field in the driven cavity example. The black circles are the results from Ghia \etal \cite{ghia1982high} and the red stars are from Erturk \etal \cite{erturk2005numerical}. The line results are taken from our solver. This simulation was run with a minimum and maximum quadtree level of 8 and 10, respectively, and a CFL of 1.0.}
        \label{fig:driven_cavity}
\end{figure}

The numerical simulations are performed using minimum and maximum quadtree levels of 8 and 10, respectively. We set the CFL number to 1.0 and dynamically refine the mesh based on the procedure discussed in Section \ref{sec:refinement}. In Figure \ref{fig:driven_cavity} we plot the results of our simulation with data from the numerical experiments of Ghia \etal \cite{ghia1982high} and Erturk \etal \cite{erturk2005numerical}. We see excellent agreement between our simulation and the previous studies.

\subsection{Oscillating cylinder} \label{sec:oscillating_cylinder}
We demonstrate the ability of our solver to handle moving interfaces by considering an oscillating cylinder example. The oscillating cylinder is a popular example of a one-way fluid-structure interaction problem that has been studied in both experimental and numerical contexts \cite{dutsch1998low, iliadis1998oscillating, seo2011sharp, GUITTET2015215}. In this test, oscillatory motion is imposed on a cylinder in a quiescent fluid and the flow around the cylinder is characterized by the Reynolds number and the Keulegan–Carpenter number defined as
\begin{align*}
    KC = \frac{u_{max}}{2 r f}.
\end{align*}

We follow the problem definition of \cite{dutsch1998low} and \cite{GUITTET2015215} and set these dimensionless parameters to $Re=100$ and $KC=5$. We set the radius of the cylinder to $r=0.05$, our frequency to $f=1$, and the amplitude of oscillation to $x_0 = 1.5914 r$. We consider the computational domain of $\Omega = [-1,1]^2$ with homogeneous Dirichlet boundary conditions and define the center position of our cylinder as
\begin{align*}
    x_c = x_0 (1 - \cos(2 \pi f t)).
\end{align*}
We compare the results of our solver with those of \cite{dutsch1998low} and \cite{GUITTET2015215} by computing the drag coefficient. We compute the drag and lift forces on the disk by integrating
\begin{equation}
\textbf{F} = \int\limits_\Gamma (-p\mathbf{I}+2 \mu \vectorr{\sigma}) \cdot\textbf{n} \label{eq:drag_integral}
\end{equation}
where $\Gamma$ is the boundary of the cylinder, $\vectorr{\sigma}$ is the symmetric stress tensor and $\textbf{n}$ is the outward facing normal vector to the cylinder. The drag and lift coefficients, $C_d$ and $C_l,$ are computed by re-scaling the components of $\vectorr{F}$ by $\rho r \norm{\vectorr{u}}_{\infty}^2.$ 
Numerically, the integral is computed using the quadrature rules of \cite{Min2007GeometricIO} on third-order smooth extensions \cite{ASLAM2004349} of the integrands.

\begin{figure}[h!]
         \includegraphics[width=0.47\textwidth]{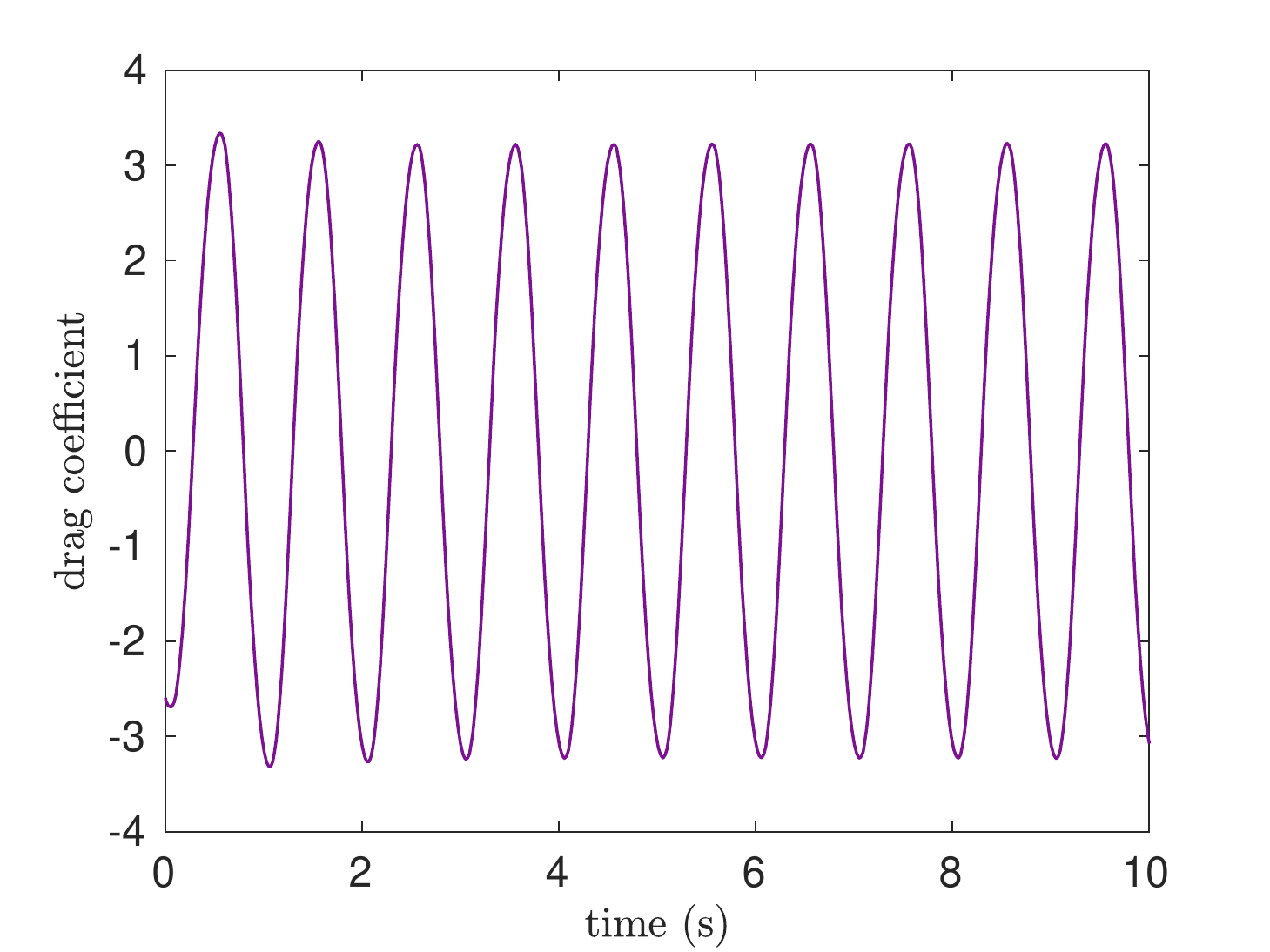}
         \includegraphics[width=0.47\textwidth]{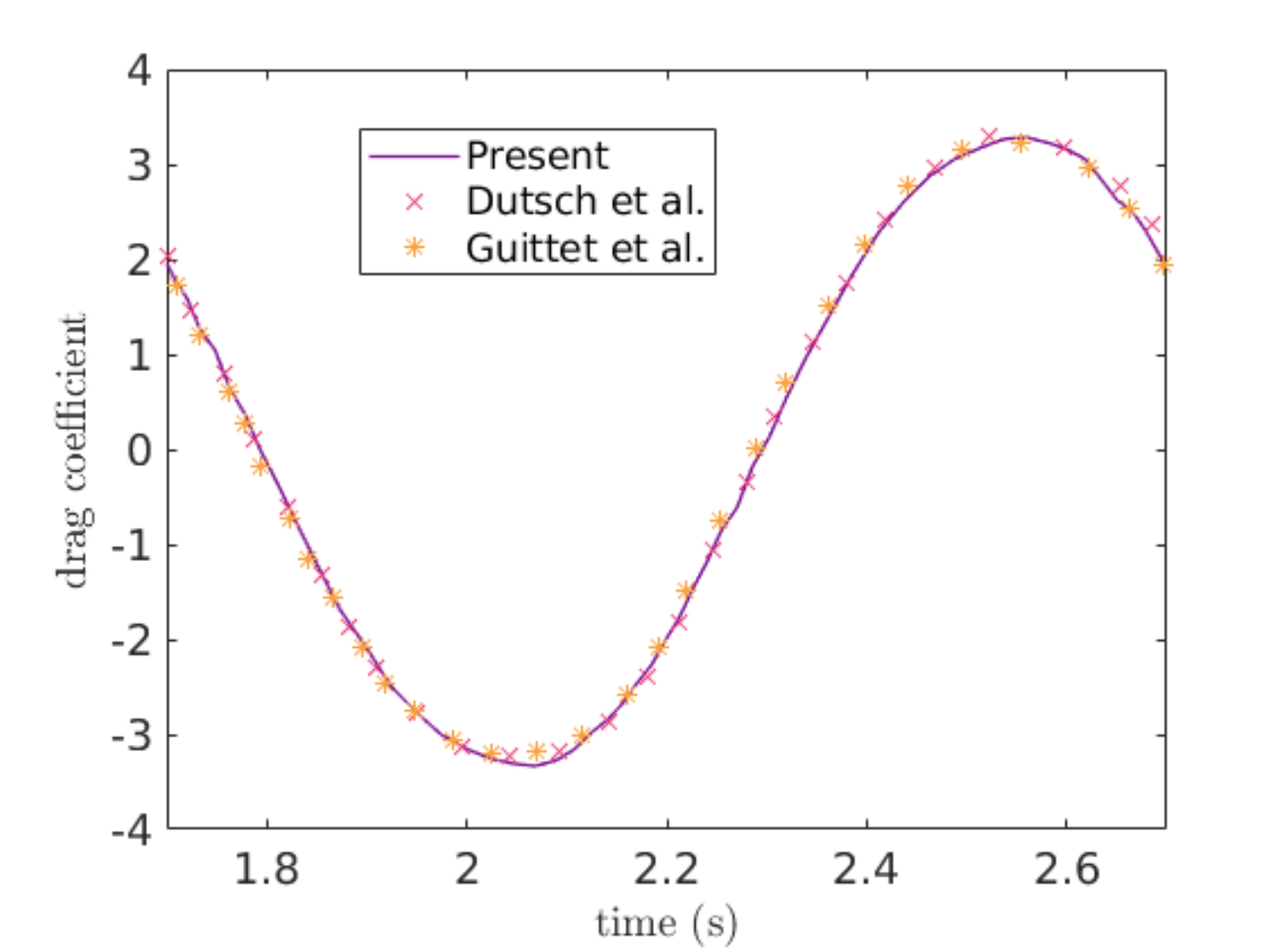}
    \caption{History of the drag coefficient for an oscillating cylinder with Re = 100 and KC = 5. The red crosses are the original experimental data from Dutsch et al. \cite{dutsch1998low} and the orange stars were computed by Guittet \etal \cite{GUITTET2015215}. Data points were digitized using the original figures. This simulation was run with a minimum and maximum quadtree level of 4 and 11, respectively, and a CFL of 10.0.}
        \label{fig:oscillating_cylinder}
\end{figure}

\begin{figure}[h!]
     \includegraphics[width=0.95\textwidth]{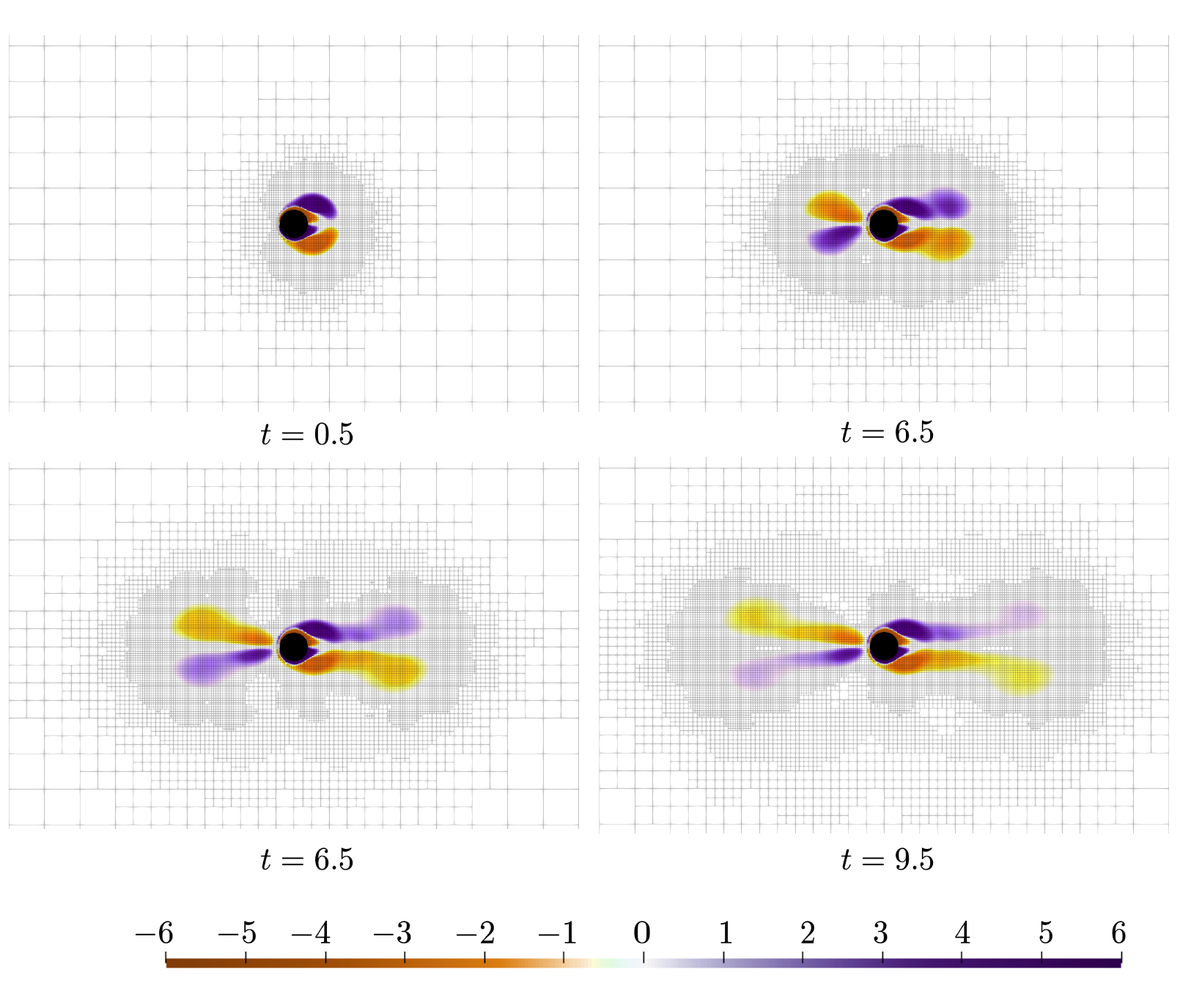}
    \caption{Time evolution of the vorticity and adaptive grid refinement for the oscillating cylinder with $Re=100$ and $KC=5$.}
        \label{fig:oscillating_cylinder_vort}
\end{figure}

For this example, we are able to use an adaptive mesh with minimum and maximum quadtree levels of 4 and 11, respectively. Additionally, we see good quantitative agreement with both studies for CFL numbers up to 10. Figure \ref{fig:oscillating_cylinder} shows the time history of the drag coefficient for $Re=100$ and $KC=5$ and a comparison with previous studies using a CFL of 10. In Figure \ref{fig:oscillating_cylinder_vort}, we show the time evolution of the vorticity and the adaptive grid refinement of the oscillating cylinder. 

\subsection{Von Kármán vortex street} \label{sec:karman_vortex}
Next, we validate our solver by considering the analysis of flow past a cylinder in two and three spatial dimensions.

\subsubsection{Two-dimensions}\label{sec:2dflowpastsphere}
The setup for this problem consists of a cylinder of radius $r = 0.5$ with its center at the location $(8,0)$ in a rectangular domain $\Omega = [0,32]\times[-8,8].$ For the velocity, we impose a Dirichlet boundary condition of $u = 1$ on the left, top, and bottom walls and a homogeneous Neumann boundary condition on the right wall. We impose homogeneous Neumann boundary conditions on the left, right, and bottom walls for the Hodge variable and a Dirichlet boundary condition of $\Phi = 0$ on the right wall. For this case, we ran our simulation on adaptive quadtree grids with minimum and maximum refinement levels of 7 and 10, respectively. To solve the advective term of the momentum equation, we use the classical departure point reconstruction method described in Section \ref{sec:improved_sl}. We also set a uniform band around the disk to properly capture the boundary layer around the interface.

\begin{figure}[h!]
    
    \centering
    \includegraphics[width = \textwidth]{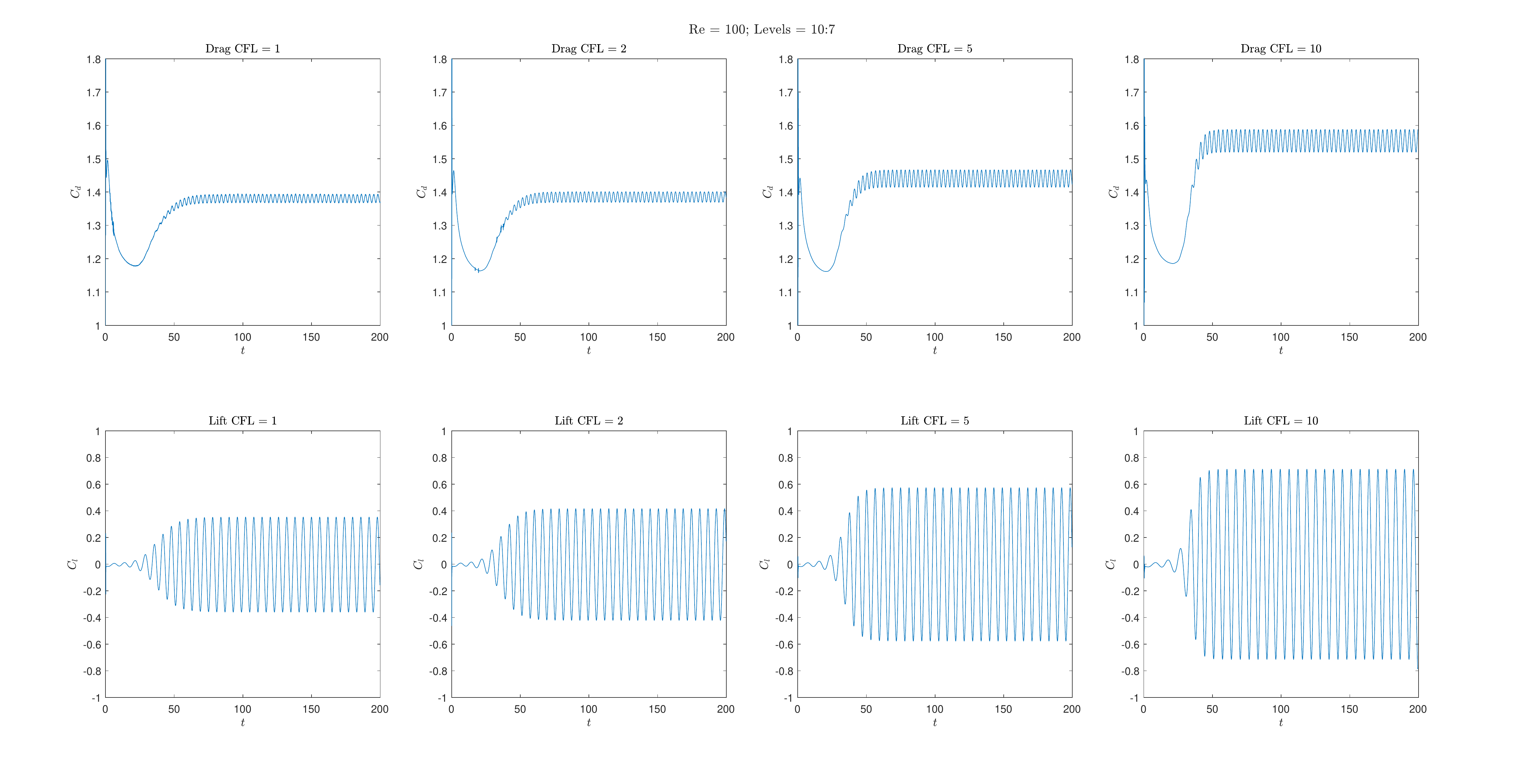} 
    \caption{Drag and lift computations for the 2-dimensional flow past disk problem with $Re = 100$ for CFL numbers of 1, 2, 5, and 10. }
    \label{fig:flowpastdisk_multiple_cfls}
\end{figure}

\begin{figure}[h!]
    \centering
   \includegraphics[width = .9\textwidth]{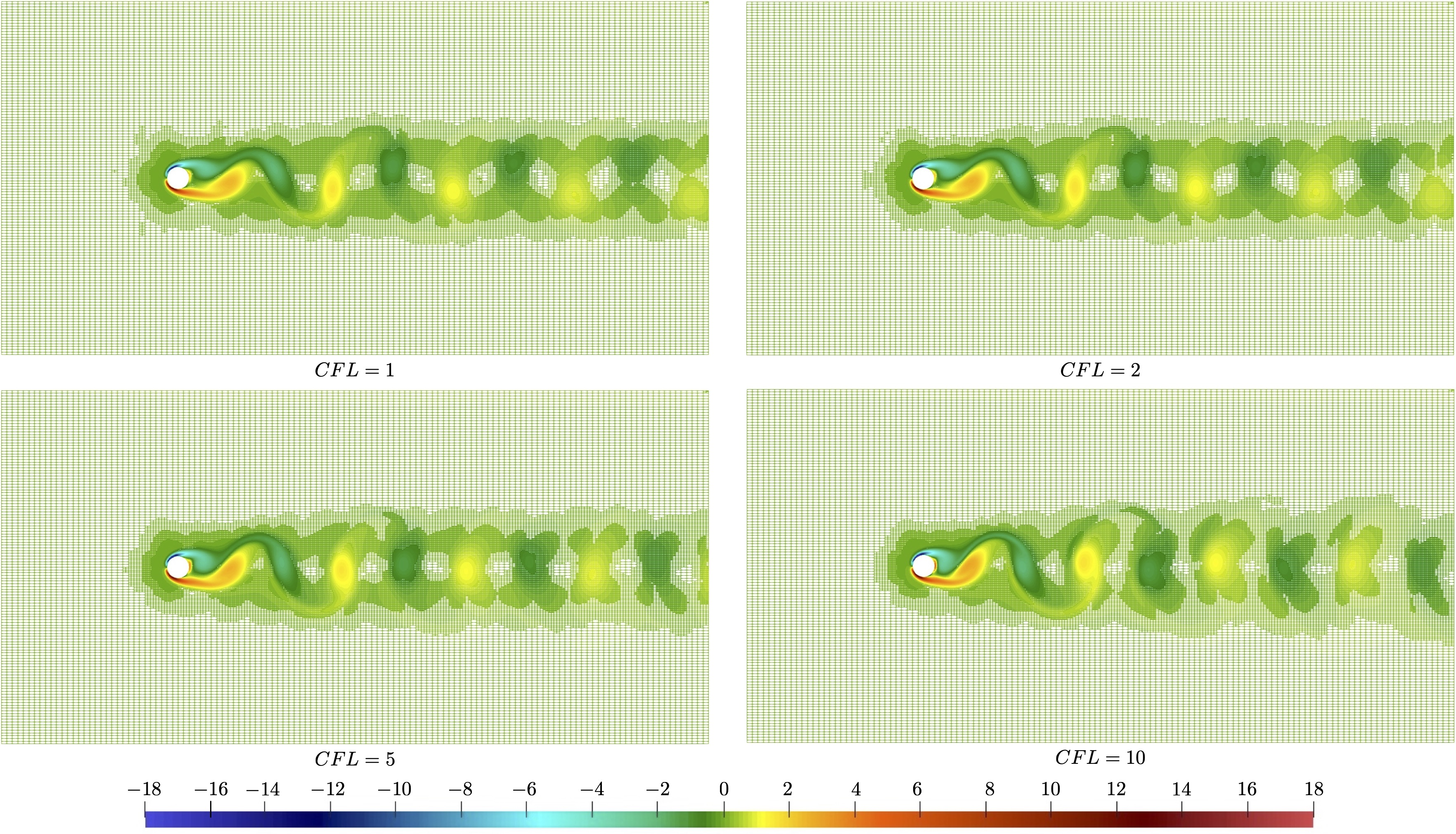}
    \caption{Visualization of the 2-dimensional flow past disk problem with \text{Re} = 100 and levels = 10:7. The vorticity and grid are  shown.  }
    \label{fig:flowpastdisk_vorticity}
\end{figure}

 We solve this problem for several CFL numbers to quantify the impact the time step has on the drag and lift forces on the cylinder and to act as a calibration tool for more complex two- and three-dimensional problems. In Figure \ref{fig:flowpastdisk_multiple_cfls}, we see that as we increase the CFL number, both the average drag coefficient values and the amplitude of the coefficient oscillations increase. Qualitatively (Figure \ref{fig:flowpastdisk_vorticity}), we observe a more severe destabilization of the wake past the cylinder for the higher CFL numbers, which is consistent with the increased forces. Nonetheless, we see that our solver can produce good qualitative results for large CFL numbers. We also see that the drag and lift coefficients are consistent with previous numerical and experimental results (Table \ref{tab:2ddragtable}).

\begin{table}[h!]
    \caption{\label{tab:2ddragtable} Drag and lift coefficients for the two-dimensional flow past disk problem.}
    \centering
    \begin{tabular}{|l|ll|ll|}\hline
                                                     & \multicolumn{2}{c|}{Drag coefficient} & \multicolumn{2}{c|}{Lift coefficient} \\
                                                     &  $Re = 100$ & $Re = 200$  & $Re = 100$ & $Re = 200$  \\ \hline
    Ng \etal  \cite{NG20098807}                     &  $1.368 \pm .016 $ & $1.373\pm0.050$ & $\pm0.360$ & $\pm0.724$\\
    Braza \etal \cite{braza_chassaing_minh_1986}     &  $1.364 \pm .015 $ & $1.400\pm 0.050$ & $\pm0.250$ & $\pm0.750$\\
    Calhoun \cite{CALHOUN2002231}                    &  $1.330 \pm .014 $ & $1.172\pm 0.058$ & $\pm0.298$ & $\pm0.668$\\
    Engelman \etal \cite{ENGELMAN1990}              &  $1.411 \pm .010 $ & -               & $\pm0.350$ &  -\\
    Guittet \etal  \cite{GUITTET2015215}            &  $1.401 \pm .011 $ & $1.383\pm 0.048$ & $\pm0.331$ & $\pm0.705$\\
    Present                                          &  $1.387 \pm .019 $ & $1.370\pm 0.060$ & $\pm0.346$ & $\pm0.762$\\
    \hline
    \end{tabular}
\end{table}

\subsubsection{Impact of departure point reconstruction}

\begin{figure}[h!]
    \centering
    \includegraphics[width = \textwidth]{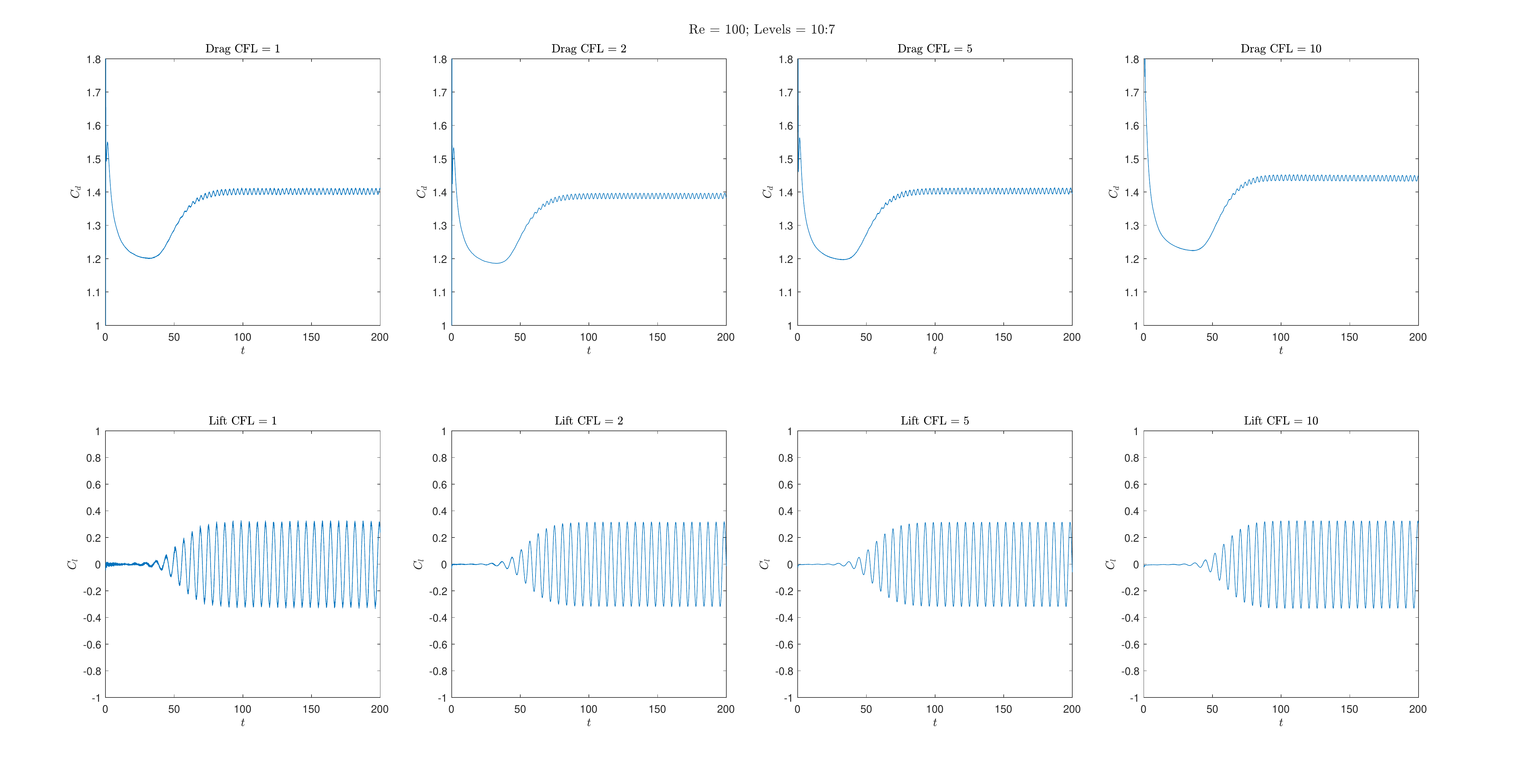} 
    \caption{Drag and lift computations for the 2-dimensional flow past disk problem with $Re = 100$ for CFL numbers of 1, 2, 5, and 10, with the improved departure point reconstruction.}
    \label{fig:flowpastdisk_multiple_cfls_updated}
\end{figure}

We implement the two-dimensional flow past disk problem in Section \ref{sec:2dflowpastsphere}, with the only difference being that we use our improved departure point method. For $Re = 100,$ we see a significant improvement in both the drag and lift coefficients for higher CFL numbers (Figure \ref{fig:flowpastdisk_multiple_cfls_updated}) compared to the implementation found in Section \ref{sec:2dflowpastsphere} (Figure \ref{fig:flowpastdisk_multiple_cfls}). 

\begin{figure}[h!]
		\centering
		\includegraphics[width = .45\textwidth]{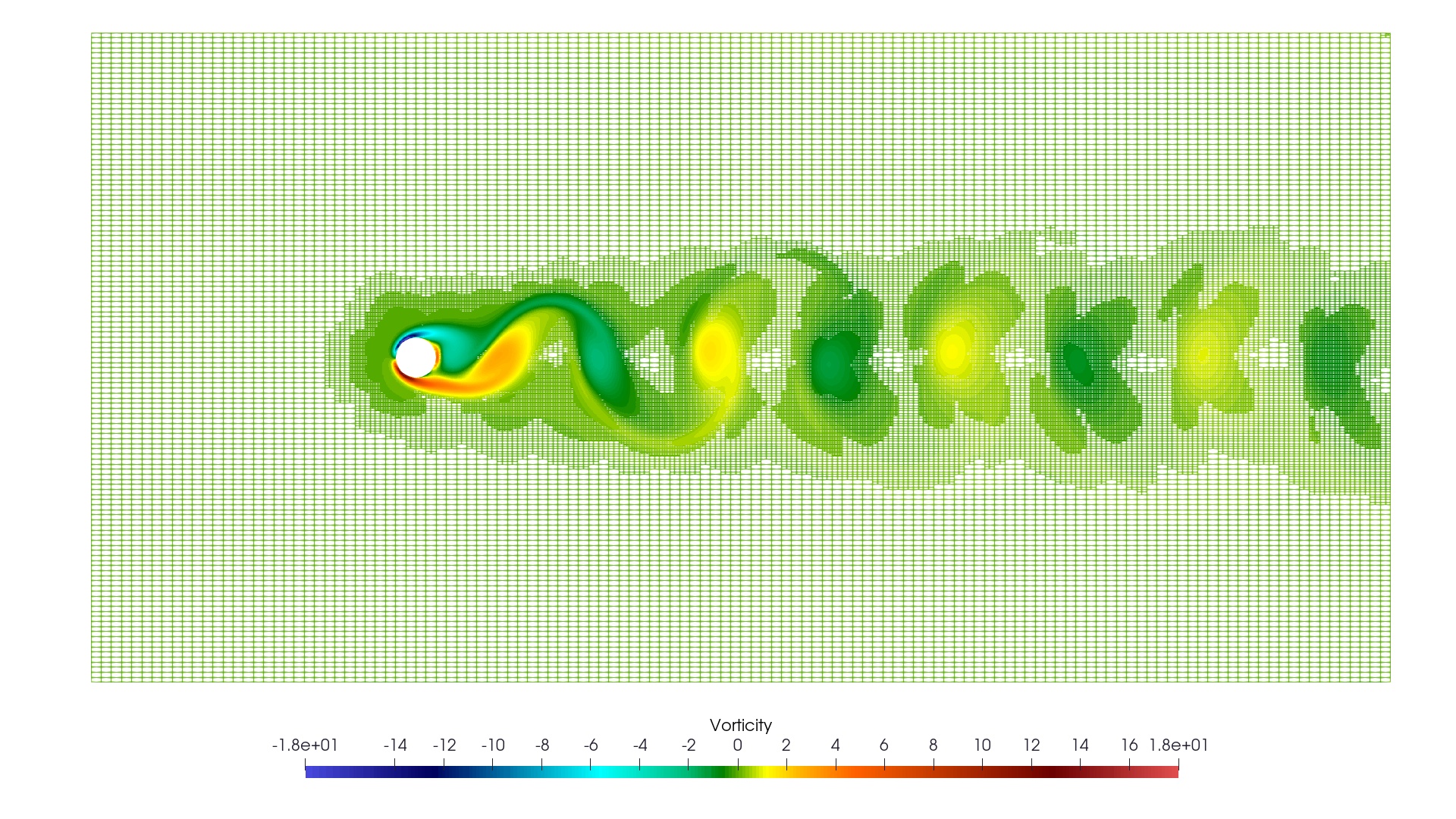}
		\includegraphics[width = .45\textwidth]{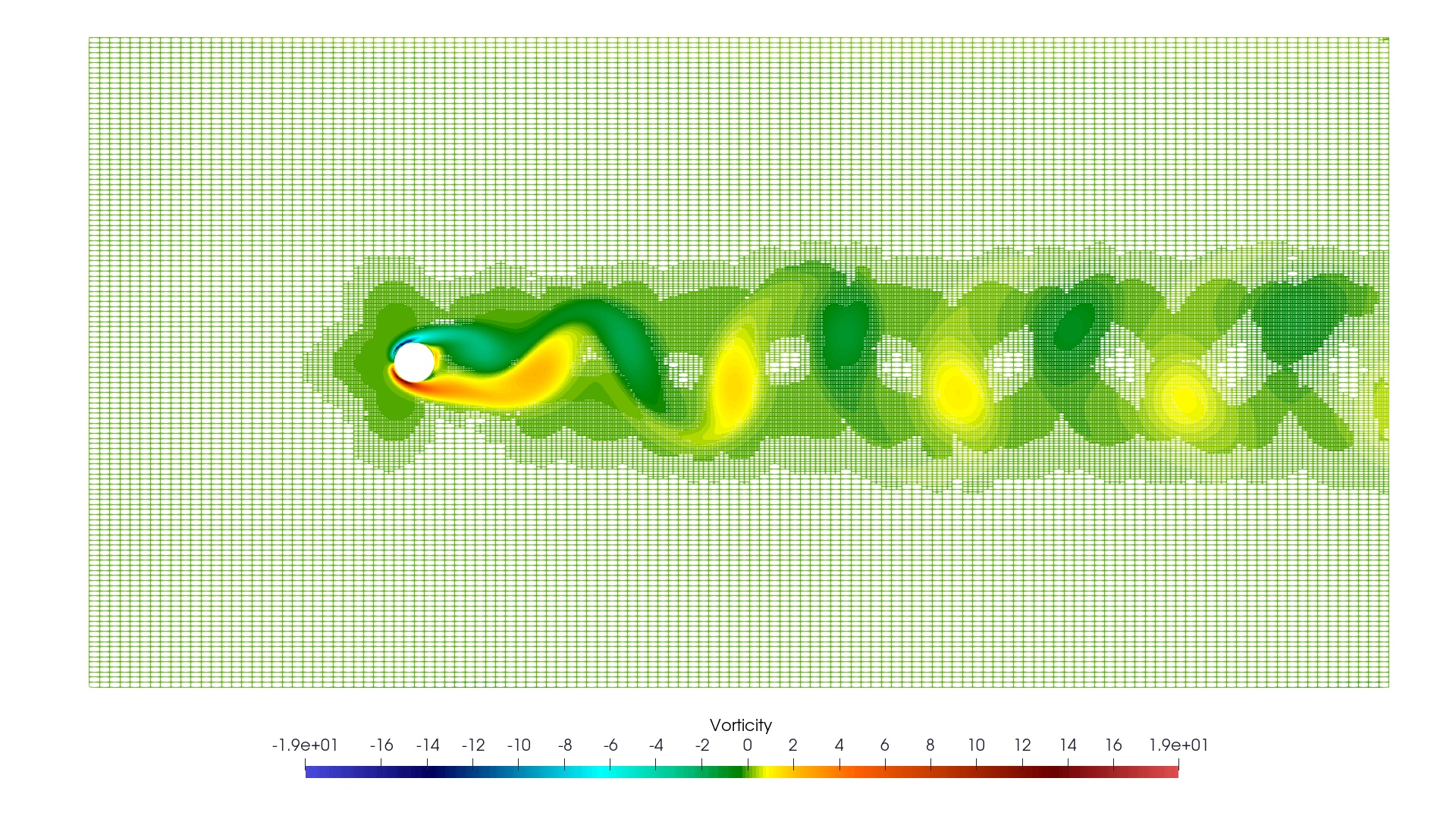}
		\caption{Visualization of the vorticity for CFL = 10. The left is the original simulation from Figure \ref{fig:flowpastdisk_vorticity}, and the right is the improved departure point reconstruction. }
		\label{fig:flowpastdisk_vorticity_updated}
	\end{figure}

These significant improvements are further realized when looking at the data qualitatively. In Figure \ref{fig:flowpastdisk_vorticity_updated}, we compare a snapshot of the vorticity profile with $CFL = 10$ of the two departure point reconstruction methods. The destabilization is not as present in the improved departure point method (as expected from the computed drag and lift computation) and looks similar to the more accurate, lower CFL number simulations.

\subsubsection{Three-dimensions}\label{sec:karman3d}
\begin{figure}[h!]
    \centering
    \includegraphics[width = .9\textwidth]{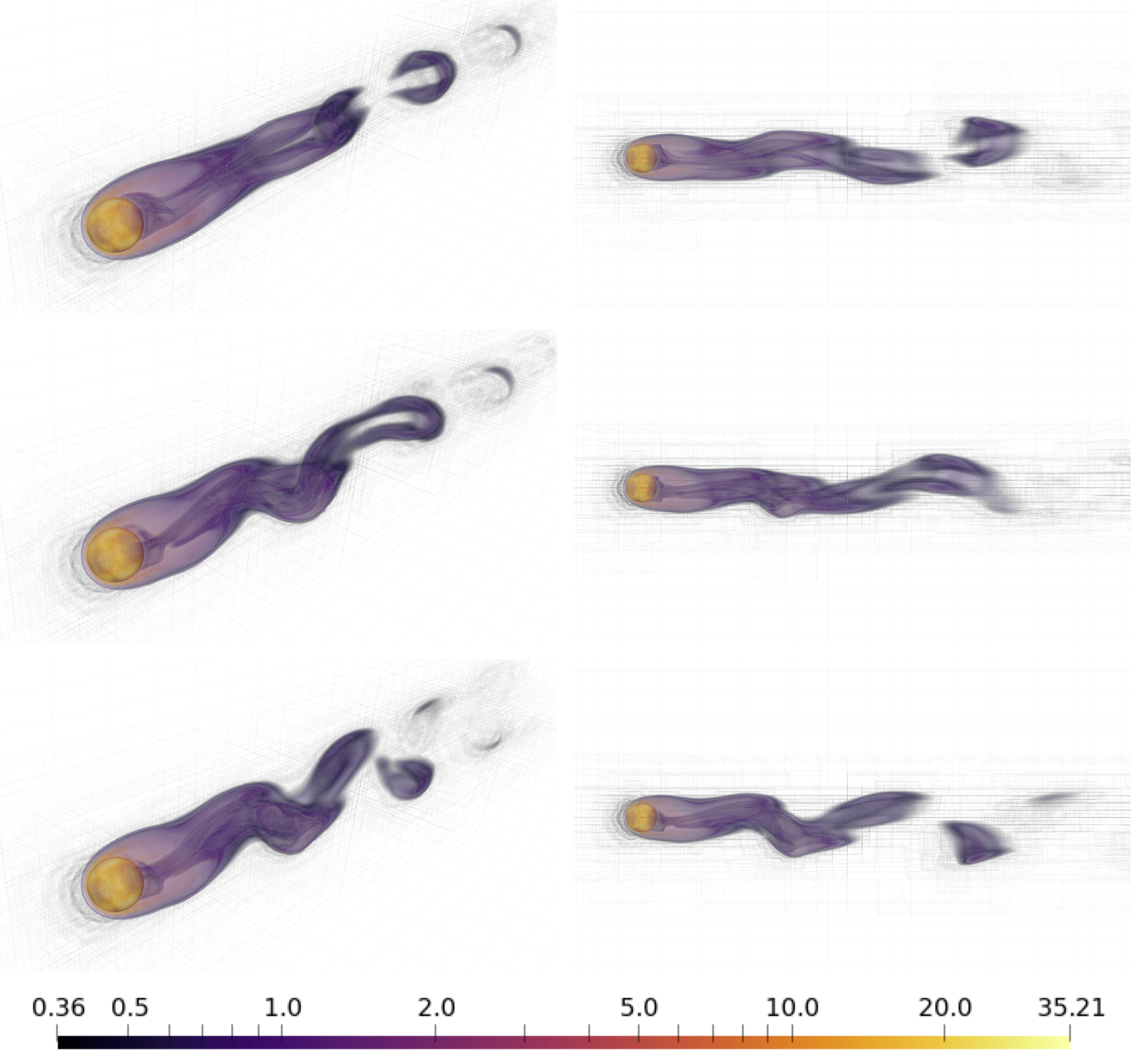}
    \caption{Visualization of flow past a sphere for $Re = 350$ at several times. The vorticity and grid are shown.}
    \label{fig:flow_past_sphere_cool}
\end{figure}

To further validate the solver and demonstrate our solver's capabilities in three dimensions, we solve the three-dimensional flow past sphere problem. The setup is similar to the 2D example in Section \ref{sec:2dflowpastsphere}, where the radius $r = 0.5$ and is centered at $(8,0,0)$ in the domain $\Omega = [0,32]\times[-8,8]\times[-8,8],$ with minimum and maximum refinement levels of 4 and 10, respectively. Like in two dimensions, we quantify the capabilities of our solver by computing the drag coefficients in a similar way to equation \eqref{eq:drag_integral}. These results can be seen in Table \ref{tab:3ddragtable}, which shows that our data agree with previous numerical and experimental studies. We also visualize the $Re = 350$ simulation in Figure \ref{fig:flow_past_sphere_cool}.

\begin{table}[!h]
\begin{center}
\caption{\label{tab:3ddragtable} Drag coefficients for the three-dimensional flow past sphere problem.}
\begin{tabular}{|l|llll|}\hline 
                                            & $Re = 100$  & $Re = 200$ & $Re  = 300$  & $Re = 350$  \\ \hline
Mittal \etal \cite{MITTAL20084825}          & 1.08        & -          & 0.67         & 0.62        \\
Marella \etal \cite{MARELLA20051}          & 1.06        & -          & 0.621        & -           \\
Le Clair \etal \cite{le1970numerical}       & 1.096       & 0.772      & 0.632        & -           \\
Johnson and Patel \cite{johnson_patel_1999} & 1.1         & 0.8        & 0.656        & -           \\
Guittet \etal  \cite{GUITTET2015215}      & 1.112       & 0.784      & 0.659        & 0.627       \\
Egan \etal \cite{EGAN2021110084}           & 1.11        & -     & 0.673        & 0.633       \\
Present                                     & 1.058       & 0.768      & 0.646        & 0.601       \\\hline
\end{tabular}\end{center}
\end{table}

\subsection{Flow past Beginnings}\label{sec:beginnings}

\begin{figure}[h!]
    \centering
    \includegraphics[width = .9\textwidth]{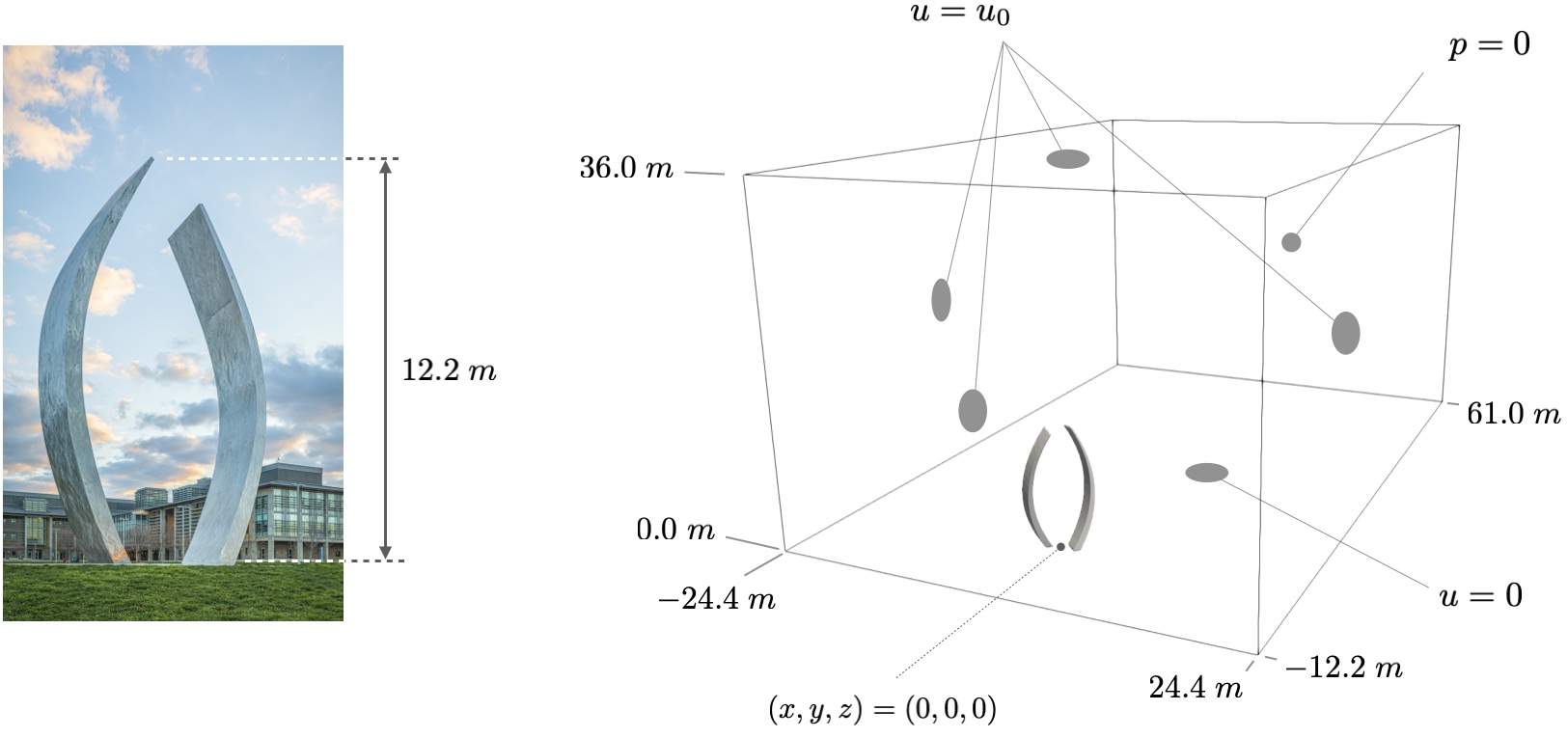}
    \caption{Problem schematic for the flow past Beginnings example with a photograph of the statue at UC Merced shown on the left.}
    \label{fig:beginnings_schematic}
\end{figure}

We conclude this section by simulating a high Reynolds number flow past the UC Merced Beginnings sculpture. The model of the sculpture was placed in the domain of $\Omega = [\SI{-12.2}{\metre},\SI{61.0}{\metre}] \times [\SI{-24.4}{\metre},\SI{24.4}{\metre}] \times [\SI{0.0}{\metre},\SI{36.8}{\metre}]$, with minimum and maximum octree refinement levels of 1 and 11, respectively. The Beginnings sculpture is approximately \SI{12.2}{\metre} (\SI{40}{ft}) in height and we center the sculpture at the point, $(x,y,z) = (0,0,0)$. We set the free-stream velocity on the boundary of the domain to $u_0=\SI{21.5}{\metre\per\second}$ (\SI{48}{mph}), which represents a windy day at UC Merced. The density of the fluid in the simulation is set to $\rho=\SI{1.3}{\kilogram\per\cubic\metre}$ and the viscosity to $\mu=\SI{1.8e-5}{\kilogram\per\metre\per\second}$, which corresponds to a Reynolds number of $Re=1.87\times10^{7}$. The velocity at the bottom wall of the domain, $z_{min}$, is set to $0$ and we set the outflow wall of the domain, $x_{max}$, to a pressure value of $0$. A schematic of the problem setup is shown in Figure \ref{fig:beginnings_schematic}. 

\begin{figure}[h!]
    \centering
    \includegraphics[width = \textwidth]{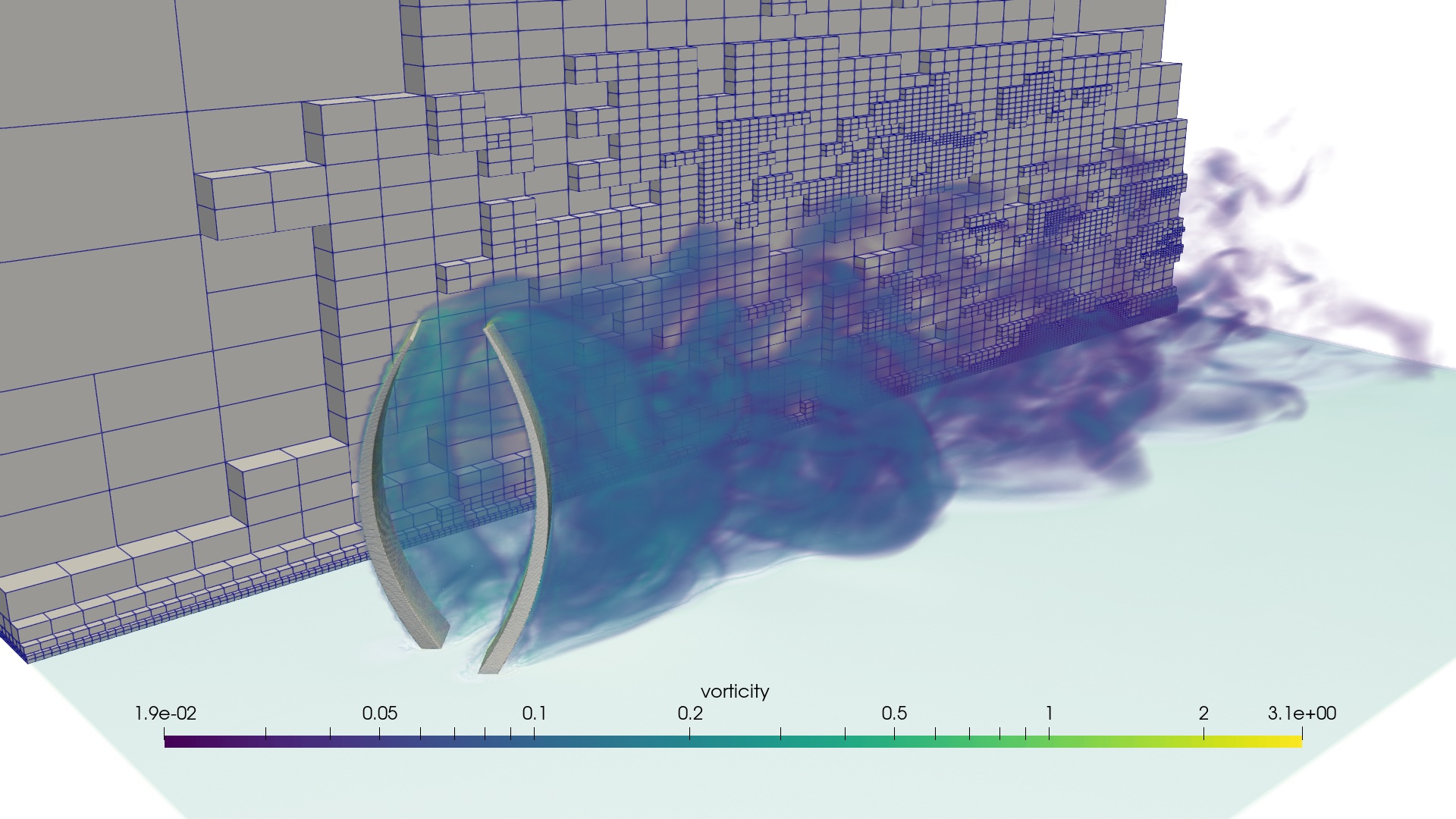}
    \caption{Visualization of the flow past the UC Merced Beginnings sculpture for $Re=1.87\times 10^{7}$. A slice of the grid is shown .}
    \label{fig:flow_past_beginnings}
\end{figure}

This example demonstrates the ability of our solver to simulate a high Reynolds number flow around an irregular geometry representing a real-world object, akin to a real work science and engineering problem. Using a 40-core machine, our nodal projection method takes approximately 560 seconds per time iteration with approximately 3~325~000 nodes using the default values for $K_{max}$ and error thresholds. Figure \ref{fig:flow_past_beginnings}, taken at $t=3~032$ seconds, demonstrates that our approach can resolve a realistic flow field around the Beginnings statue, including the boundary layer while maintaining stability and computational efficiency. 

\section{Conclusions} \label{sec:ccl}
We presented a novel projection method for the simulation of incompressible flows in arbitrary domains using quad/octrees, where all of the variables are collocated at the grid nodes. 
By design, our collocated projection operator is an iterative procedure. If it exists, the limit of this iterated procedure is the canonical projection on the incompressible space (Section \ref{sec:pres_and_gen_prop}).
Leveraging the connection between the collocated and MAC layouts, we found a technical sufficient condition for this limit to exist on periodic grids. 
We then verified the stability and convergence of our projection operator numerically in the presence of various boundary and interface conditions. Next, we validated our nodal Navier-Stokes solver using canonical two- and three-dimensional examples and see that our collocated solver achieves high-order convergence (Section \ref{sec:avortex}), replicates experimental data (Section \ref{sec:oscillating_cylinder}), and can accurately represent standard computational results with CFL numbers far greater than 1 (Section \ref{sec:karman_vortex}). Finally, we showed that our solver is well suited to study real-world science and engineering problems by simulating high Reynolds number flows past arbitrary geometries (Section \ref{sec:beginnings}). Our nodal projection method has proven to be a competitive computational fluid dynamics tool for studying complex fluid flows. 

However, we see room for further exploration and improvement at both theoretical and computational levels. Although we derived a technical sufficient condition for the stability of our projection operator, we could not formally prove its stability for our specific ghost node construction and left the boundary effects out of our analysis. Therefore, we relied on an in-depth computational verification process. A more detailed theoretical study could provide information to optimize the construction of collocated operators. Furthermore, our implementation was parallelized using a shared-memory framework. To further extend the capabilities of our nodal solver, we aim to expand it to a distributed memory framework, as this significant enhancement will provide stronger scalability and enable us to investigate a broader range of applications. 

To conclude, we emphasize that our nodal projection method uses a fully collocated data layout, drastically reducing implementation costs. The results presented herein demonstrate that our method is a suitable foundation to study a wide range of fluid flows, and we believe that this work can be readily extended to solve more complex phenomena, such as free-surface and multi-phase flows. We foresee our nodal projection method leading to the development of a new class of collocated computational tools capable of accelerating scientific discoveries by lowering accessibility barriers.

\section{Acknowledgments}
The authors thank Brandon Stark for collecting and processing the data to generate a three-dimensional model of the UC Merced Beginnings statue. This material is based upon work supported by the National Science Foundation under Grant No. DMS-1840265.

\section*{Computational Environment}
Our solver was implemented in C++ 20, and all computations were done in parallel on CPUs using OpenMP \cite{OpenMP_paper,openmp18} libraries. Three-dimensional visualizations were generated using ParaView \cite{Ahrens2005ParaViewAE} (version 5.11.0). The imagery of the Beginnings sculpture was collected using a Skydio 2+ drone and processed with DroneDeploy using default settings to generate a three-dimensional model.

\section*{Credit author statement }  
{\bf Matthew Blomquist}: 
Formal analysis, Investigation, Software, Validation, Visualization, Writing - Original Draft, Writing - Review $\&$ Editing. 
 {\bf Scott R. West}: 
Formal analysis, Investigation, Software, Validation, Visualization, Writing - Original Draft, Writing - Review $\&$ Editing.
{\bf Adam L. Binswanger}: 
Formal analysis, Investigation, Software, Validation, Writing - Original Draft, Writing - Review $\&$ Editing.
{\bf Maxime Theillard}: Conceptualization, Supervision, Formal Analysis, Methodology, Software, Project administration, Writing - Original Draft, Writing - Review $\&$ Editing. 
 
\bibliographystyle{abbrv}
\bibliography{main.bib}

\end{document}